\documentclass[11pt,a4paper]{article}
\usepackage[T1]{fontenc}
\usepackage[utf8]{inputenc}
\usepackage{lmodern,amsmath,amssymb,amsthm,mathtools,microtype}
\usepackage[margin=28mm]{geometry}
\usepackage{xurl,pgfplots}
\pgfplotsset{compat=1.18}
\usepackage[authoryear,round,longnamesfirst]{natbib}
\usepackage[hidelinks,hyperfootnotes=false]{hyperref}
\hypersetup{pdftitle={The Prodigal Elephant Returns: Last-Return Laws in the Superdiffusive Regime},pdfauthor={Hélio Trinas}}
\numberwithin{equation}{section}
\newtheorem{theorem}{Theorem}[section]
\newtheorem{proposition}[theorem]{Proposition}
\newtheorem{lemma}[theorem]{Lemma}
\newtheorem{corollary}[theorem]{Corollary}
\theoremstyle{definition}

\theoremstyle{remark}

\newcommand{\Pp}{\mathbb P}
\newcommand{\E}{\mathbb E}
\newcommand{\R}{\mathbb R}
\newcommand{\Z}{\mathbb Z}
\newcommand{\F}{\mathcal F}
\newcommand{\LL}{\mathcal L}
\newcommand{\ind}{\mathbf 1}
\newcommand{\dd}{\,\mathrm d}
\newcommand{\law}{\xrightarrow{\mathrm d}}
\newcommand{\inprob}{\xrightarrow{\mathrm p}}
\newcommand{\asconv}{\xrightarrow{\mathrm{a.s.}}}
\DeclareMathOperator{\Var}{Var}
\DeclareMathOperator{\sgn}{sgn}
\allowdisplaybreaks[2]
\title{The Prodigal Elephant Returns:\\Last-Return Laws in the Superdiffusive Regime}
\author{H\'elio Trinas\thanks{E-mail: helio@dme.ufrj.br. Supported by CAPES.}}
\date{}
\begin{document}
\maketitle
\begin{abstract}
Let $(S_n)$ be the one-dimensional elephant random walk with memory parameter $p\in(3/4,1)$ and first-step probability $q=\Pp(S_1=1)\in[0,1]$. Set $a=2p-1$, and let $f_{q,a}$ denote the density of the almost-sure limit of $S_n/n^a$. We prove that the last return time $T_0=\sup\{m\geq0:S_m=0\}$ satisfies
\[
 \Pp(T_0\geq n)\sim f_{q,a}(0)\sqrt{\frac{8}{\pi(2a-1)}}\,n^{-(a-1/2)}.
\]
Reflection gives $f_{q,a}(0)=f_{1,a}(0)$, so the constant is independent of $q$. For future walks started from admissible states $(n,x_n)$ with $x_n/\sqrt n\to z$, we establish functional convergence of the centered martingale-normalized paths to a time-changed Brownian motion. Uniform control of late hitting times then gives the eventual-return kernel $2\Phi(-\sqrt{2a-1}\,|z|)$, uniformly on bounded scaled windows, where $\Phi$ is the standard normal distribution function. The same last-visit asymptotic holds uniformly for nonempty deterministic targets $F_n\subseteq[-r_n,r_n]\cap\Z$, $r_n=o(\sqrt n)$, held fixed during each trajectory. Conditionally on $T_0\geq n$, we identify the joint limit of $S_n/\sqrt n$, the next return divided by $n$, and the last return divided by the next return. The last ratio converges to $A^{-1/(2a-1)}$, independently of the other two limiting coordinates, where $A$ has the arcsine distribution on $(0,1)$.
\end{abstract}
\begingroup
\renewcommand{\thefootnote}{}
\endgroup

\section{Introduction}
The Elephant Random Walk (ERW), introduced by \citet{ST}, is a nearest-neighbor walk on $\Z$, started from $S_0=0$, whose next increment repeats a uniformly selected past increment with probability $p$ and reverses it otherwise. The first increment is positive with probability $q\in[0,1]$. The limit theorems of \citet[Theorems 1--3]{BB} and \citet[Theorems 2--3]{CGS} distinguish three scaling regimes: the diffusive scale $\sqrt n$ for $p<3/4$, the critical scale $\sqrt{n\log n}$ for $p=3/4$, and the superdiffusive scale $n^{2p-1}$ for $3/4<p<1$. \citet{BB} use the connection with P\'olya-type urns; martingale methods appear in \citet[Section III.A]{CGS} and \citet{Bercu}.

The walk is recurrent for $p\leq3/4$ and transient for $p>3/4$; see \citet{CP} and \citet[Theorem 1.3]{Qin}. Return times raise finer questions than this classification. \citet{Bertoin} studies zeros and first returns in the diffusive regime, while \citet{FangCritical} treats the critical regime, including the last zero before a finite observation horizon. We consider the last zero over the entire trajectory in the transient regime: how late can it occur, and what does the walk look like when it still has a return ahead of it?

Throughout the paper,
\[
 3/4<p<1,\qquad a=2p-1\in(1/2,1),\qquad \delta=2a-1.
\]
The normalized position $S_n/n^a$ converges almost surely to a random variable $\LL_{q,a}$; see \citet[Theorem 3]{BB}, \citet[Theorem 3]{CGS}, and \citet[Theorem 3.7]{Bercu}. \citet[Theorem 1.3]{GLR} prove that its law has a positive bounded smooth density $f_{q,a}$, while \citet{GLRS} further investigate this limit law and its moments. We keep $q$ arbitrary throughout and write $f_a=f_{1,a}$ for the density corresponding to a positive first step. Reflection gives $f_{q,a}(0)=f_a(0)>0$ for every $q\in[0,1]$, so the density factor in the last-return constant is independent of the first-step bias. Almost-sure convergence and absolute continuity together yield a direct proof of the known transience result \citep{CP,Qin}, recalled in Section \ref{sec:model}.

Late returns are mainly contributed by trajectories whose position at time $n$ is comparable to the fluctuations still available in the future. After removing the drift amplification, these fluctuations have variance of order $n$ in the relevant window. The return scale is therefore $\sqrt n$, smaller than the typical displacement $n^a$. The local limit theorem (LLT) of \citet[Theorem 3.1]{TV} assigns mass of order $n^{-a}$ to each of the order $\sqrt n$ sites in this window, suggesting the exponent $a-1/2$. The return probability from each scaled position supplies the weight that determines the constant.

\begin{theorem}\label{thm:tail}
Let $T_0=\sup\{m\geq0:S_m=0\}$. Then $T_0<\infty$ almost surely and
\begin{equation}\label{eq:tail}
 \Pp(T_0\geq n)\sim f_a(0)\sqrt{\frac{8}{\pi(2a-1)}}\,n^{-(a-1/2)}.
\end{equation}
The law of $T_0$ is independent of $q$. For every fixed nonempty finite set $F\subset\Z$, let $T_F=\sup\{m\geq0:S_m\in F\}$, with $\sup\varnothing=0$. Then $T_F<\infty$ almost surely and \eqref{eq:tail} holds with $T_F$ in place of $T_0$, with the same constant. For $r>0$, $\E[T_F^r]<\infty$ if and only if $r<a-1/2$.
\end{theorem}

The law of $T_F$ may depend on $q$, although its leading tail constant does not. A fixed finite target shrinks to a point on the $\sqrt n$ scale, which explains why the constant contains no factor $|F|$. Corollary~\ref{cor:subdiffusive-targets} extends the asymptotic uniformly to nonempty deterministic targets contained in $[-r_n,r_n]$ when $r_n=o(\sqrt n)$. The target may vary with the observation time $n$ but remains fixed along each trajectory.

Our second result describes a late return jointly with the position of the walk. Define
\[
 \tau_n=\inf\{m\geq n:S_m=0\},\qquad E_n=\{T_0\geq n\}=\{\tau_n<\infty\},\qquad n\geq1.
\]
Thus $\tau_n$ is the next zero at or after $n$. On $E_n$, the ratios $\tau_n/n$ and $T_0/\tau_n$ describe, respectively, when this next return occurs and how much later the final return occurs. Let $U,R,\varepsilon,A$ be independent, where $U$ is uniform on $(0,1)$, $R$ has the Rayleigh density $r e^{-r^2/2}\ind_{\{r>0\}}$, $\varepsilon$ takes the values $-1,+1$ with equal probabilities, and $A$ has the arcsine density $\{\pi\sqrt{u(1-u)}\}^{-1}$ on $(0,1)$.

\begin{theorem}\label{thm:joint-intro}
Conditionally on $E_n$,
\begin{equation}\label{eq:joint-intro}
 \left(\frac{S_n}{\sqrt n},\frac{\tau_n}{n},\frac{T_0}{\tau_n}\right)
 \law
 \left(\frac{\varepsilon U R}{\sqrt\delta},
 (1-U^2)^{-1/\delta}, A^{-1/\delta}\right).
\end{equation}
In particular, the third limiting coordinate is independent of the first two. The limiting distribution function of $\tau_n/T_0$, given $E_n$, is
\[
 r\longmapsto\frac2\pi\arcsin(r^{\delta/2}),\qquad 0\leq r\leq1.
\]
\end{theorem}

The shared variable $U$ encodes the dependence between the limiting position and next return. The independent arcsine factor describes the future after that return, with absolute times divided by the return time. Its classical antecedent is the arcsine law of \citet{Levy} for the last Brownian zero; see also \citet[Theorem~5.26]{MP}. For the critical ERW, \citet[Proposition 3.1(ii)]{FangCritical} obtains a Brownian last-zero limit before a finite horizon, on a logarithmic scale. The ratio in \eqref{eq:joint-intro} concerns the global last return after a late next zero in the superdiffusive regime. To the best of our knowledge, the exact tail \eqref{eq:tail} and the joint conditional law \eqref{eq:joint-intro} have not appeared previously. First-return problems with prescribed histories \citep{FangTraining} and conditioning to avoid zero \citep{IT} have been studied in the diffusive and critical regimes.

As a direct consequence of the LLT of \citet[Theorem 3.1]{TV} and the known transience \citep{CP,Qin}, we also record
\[
 \ell_\infty(x):=\sum_{m=0}^{\infty}\ind_{\{S_m=x\}}<\infty\quad\text{a.s.},
 \qquad \E[\ell_\infty(x)]=\infty,
 \qquad x\in\Z.
\]
Corollary \ref{thm:localtime} gives the growth of the mean occupation up to a finite time. This measures the number of visits; $\E[T_0]=\infty$ instead concerns the time of the last one. For the ERW, a return to zero at time $m$ does not restore the law started at time zero, since the transitions retain their dependence on absolute time. Successive excursions therefore do not form identically distributed renewal blocks, a distinction also central to Bertoin's analysis of zeros in the diffusive regime \citep[Introduction and Section 5]{Bertoin}.

The proof of the tail starts from the decomposition
\[
 \Pp(E_n)=\sum_x\Pp(S_n=x)H_n(x),
\]
where $H_n(x)$ is the probability of a future zero for the walk started from $x$ at absolute time $n$. We use the martingale normalization of \citet{CGS} and \citet{Bercu}. \citet[Sections 2 and 4]{Bertoin} studies diffusive return times through a Brownian embedding of this normalized walk, and \citet[Sections 2--3]{FangTraining} combines martingale functional limits with Brownian hitting arguments for growing prescribed histories. 

Gaussian fluctuations around the superdiffusive limit were studied by \citet{KT}; \citet[Theorem 1.1]{Tok} obtains a functional limit for future fluctuations under the original law. Here we need a functional limit for walks started at time $n$ from states $x_n/\sqrt n\to z$, which we prove in Theorem~\ref{thm:future-fclt}. We verify the martingale functional central limit theorem in the form stated by \citet[Theorem 2.1(ii)]{Whitt} for the changing initial-state laws needed here. Two further estimates are essential: one controls hits beyond a large rescaled time, and another controls starting points outside bounded $\sqrt n$ windows. Together with the LLT of \citet{TV}, they justify the above averaging and determine the constant.

Section \ref{sec:model} records the transition rule and the known limit results. Sections \ref{sec:martingale}--\ref{sec:tail} prove the hitting estimates and the tail. Sections \ref{sec:conditional}--\ref{sec:last} identify the joint law and the final escape direction. Sections \ref{sec:finite} and \ref{sec:occupation} treat finite targets, including targets growing below the $\sqrt n$ scale, and mean occupation, respectively. Appendix~\ref{app:estimates} collects the product and second-moment estimates and the verification of the explicit conditional representation.

\section{The model and known limit results}\label{sec:model}
We write $\law$, $\inprob$ and $\asconv$ for convergence in distribution, in probability and almost surely, respectively. Unlabelled arrows denote convergence of deterministic quantities; $\stackrel{\mathrm d}=$ denotes equality in distribution.
Let $X_n\in\{-1,+1\}$, with $\Pp(X_1=1)=q$, and set
\[
 S_0=0,\qquad S_n=\sum_{j=1}^nX_j,\qquad \F_n=\sigma(X_1,\ldots,X_n).
\]
Following the original definition in \citet{ST}, given $\F_n$, choose an index uniformly from $\{1,\ldots,n\}$, independently of the coin deciding whether the selected increment is repeated or reversed. Since the numbers of positive and negative increments are $(n+S_n)/2$ and $(n-S_n)/2$, respectively, the transition and drift identities used in \citet[Section 2]{Bercu} are
\begin{align}
 \Pp(X_{n+1}=1\mid\F_n)
 &=p\frac{n+S_n}{2n}+(1-p)\frac{n-S_n}{2n}
 =\frac12\left(1+a\frac{S_n}{n}\right),\label{eq:transition}\\
 \E[X_{n+1}\mid\F_n]&=aS_n/n.\label{eq:drift}
\end{align}
Thus $(S_n)$ is a time-inhomogeneous Markov chain. A state $(n,x)$ is admissible if $n\geq1$, $x\in\Z$, $|x|\leq n$ and $x\equiv n\pmod2$. For every such state, $\Pp_{n,x}$ is the future law generated by \eqref{eq:transition}, initialized at $x$ at absolute time $n$; $\E_{n,x}$ denotes its expectation. Whenever $\Pp(S_n=x)>0$, it agrees with the conditional future law given $S_n=x$. This definition is independent of $q$ and remains valid when $\Pp(S_n=x)=0$. Under $\Pp_{n,x}$, $\F_m$ denotes the natural filtration $\sigma(S_n,\ldots,S_m)$ of this future walk, for $m\geq n$. In sums, kernels are set to zero for $|x|>n$.

The role of absolute time can already be seen in a two-step excursion. For even $m\geq2$, \eqref{eq:transition} gives
\[
 \Pp_{m,0}(S_{m+2}=0)=\frac12\left(1-\frac a{m+1}\right),
 \qquad \Pp(S_2=0)=\frac{1-a}{2}=1-p.
\]
Returning to zero makes the next step symmetric, but it does not erase the time dependence of subsequent transitions.

We will use the following convergence, density and local limit results.

\begin{proposition}\label{thm:llt}
There is a random variable $\LL_{q,a}$ such that $S_n/n^a\asconv\LL_{q,a}$, by \citet[Theorem 3.7]{Bercu}. Its law has a positive bounded smooth density $f_{q,a}$ on $\R$, by \citet[Theorem 1.3]{GLR}. The uniform local limit theorem of \citet[Theorem 3.1]{TV} states that
\begin{equation}\label{eq:llt}
 \sup_{\substack{x\in\Z\\x\equiv n\ (\mathrm{mod}\ 2)}}
 \left|\frac{n^a}{2}\Pp(S_n=x)-f_{q,a}(x/n^a)\right|\longrightarrow0.
\end{equation}
Consequently, for a finite constant $C=C(a,q)$,
\begin{equation}\label{eq:anti}
 \sup_x\Pp(S_n=x)\leq Cn^{-a},\qquad n\geq1.
\end{equation}
\end{proposition}

The bound \eqref{eq:anti} follows from \eqref{eq:llt} and boundedness of $f_{q,a}$.

Conditioning on $X_1$ and reflecting all increments gives
\[
 f_{q,a}(y)=qf_a(y)+(1-q)f_a(-y),\qquad f_{q,a}(0)=f_a(0)>0.
\] The same reflection preserves the entire zero set, so its law, including that of $T_0$, is independent of $q$. Finally, on $\{\LL_{q,a}\ne0\}$ there is a random $n_0$ such that
\[
 |S_n|\geq\tfrac12|\LL_{q,a}|n^a\quad(n\geq n_0).
\]
Hence $|S_n|\asconv\infty$. Every fixed finite set is therefore visited only finitely often. This recovers the superdiffusive part of the recurrence--transience classification in \citet{CP} and \citet{Qin} from the convergence and density results cited in Proposition \ref{thm:llt}.

\section{A martingale for the future walk}\label{sec:martingale}
We use the multiplicative martingale normalization appearing in \citet[Section III.A]{CGS} and \citet[Section 2]{Bercu}, here normalized at the starting time $n$. For $m\geq n\geq1$, define
\[
 R_{n,m}=\prod_{j=n}^{m-1}\left(1+\frac aj\right)^{-1},\qquad R_{n,n}=1.
\]
This factor removes the amplification of the conditional mean between times $n$ and $m$. Multiplying by $R_{n,m}$ isolates the future fluctuations while preserving every visit to zero.
Indeed, if $b_m=\Gamma(m)\Gamma(1+a)/\Gamma(m+a)$ is the factor in \citet[(2.4)]{Bercu}, then $R_{n,m}=b_m/b_n$. Constants denoted by $c_a,C_a$ may change from one occurrence to another and depend only on $a$, unless further dependence is stated.

\begin{lemma}\label{lem:normalization}
There are $0<c_a\leq C_a<\infty$ such that
\[
 c_a(n/m)^a\leq R_{n,m}\leq C_a(n/m)^a.
\]
For every fixed $0\leq T<\infty$,
\begin{equation}\label{eq:Rlimit}
 \sup_{0\leq t\leq T}|R_{n,n+\lfloor nt\rfloor}-(1+t)^{-a}|\longrightarrow0.
\end{equation}
For $T\geq1$,
\begin{equation}\label{eq:variance-tail}
 \frac1n\sum_{m\geq n+\lfloor nT\rfloor}R_{n,m+1}^2
 \leq C_a(1+T)^{1-2a}.
\end{equation}
\end{lemma}
The product estimates are proved in Appendix~\ref{app:estimates}.

For a walk started at $(n,x)$, the martingale increment and conditional-variance formulas of \citet[Appendix A, (A.1)--(A.4)]{Bercu} take the following form.

\begin{lemma}\label{lem:martingale}
Under $\Pp_{n,x}$, the process $M_m^{(n)}=R_{n,m}S_m$, $m\geq n$, is a martingale and
\begin{align}
 \Delta M_{m+1}^{(n)}&=R_{n,m+1}(X_{m+1}-aS_m/m),\label{eq:increment}\\
 |\Delta M_{m+1}^{(n)}|&\leq2R_{n,m+1},\label{eq:jumps}\\
 \E_{n,x}[(\Delta M_{m+1}^{(n)})^2\mid\F_m]
 &=R_{n,m+1}^2\left(1-a^2S_m^2/m^2\right).\label{eq:incrementvar}
\end{align}
\end{lemma}
\begin{proof}
The drift identity \eqref{eq:drift} and $R_{n,m+1}(1+a/m)=R_{n,m}$ give the martingale property and \eqref{eq:increment}. Since $|S_m|\leq m$ and $X_{m+1}^2=1$,
\[
 |X_{m+1}-aS_m/m|\leq1+a\leq2,
 \qquad \Var(X_{m+1}\mid\F_m)=1-a^2S_m^2/m^2.
\]
These yield \eqref{eq:jumps} and \eqref{eq:incrementvar}.
\end{proof}

\begin{lemma}\label{lem:second}
For $0<K,T<\infty$ there is $C=C(a,K,T)<\infty$ such that, for every integer $n\geq1$,
\[
 \sup_{\substack{x\in\Z:\ (n,x)\text{ admissible}\\|x|\leq K\sqrt n}}
 \sup_{n\leq m\leq n+\lfloor nT\rfloor}\E_{n,x}[S_m^2]\leq Cn.
\]
\end{lemma}
The proof in Appendix~\ref{app:estimates} uses the second-moment recursion of \citet[(14)]{ST}; its conditional form is given in \citet[(A.3)]{Bercu}.

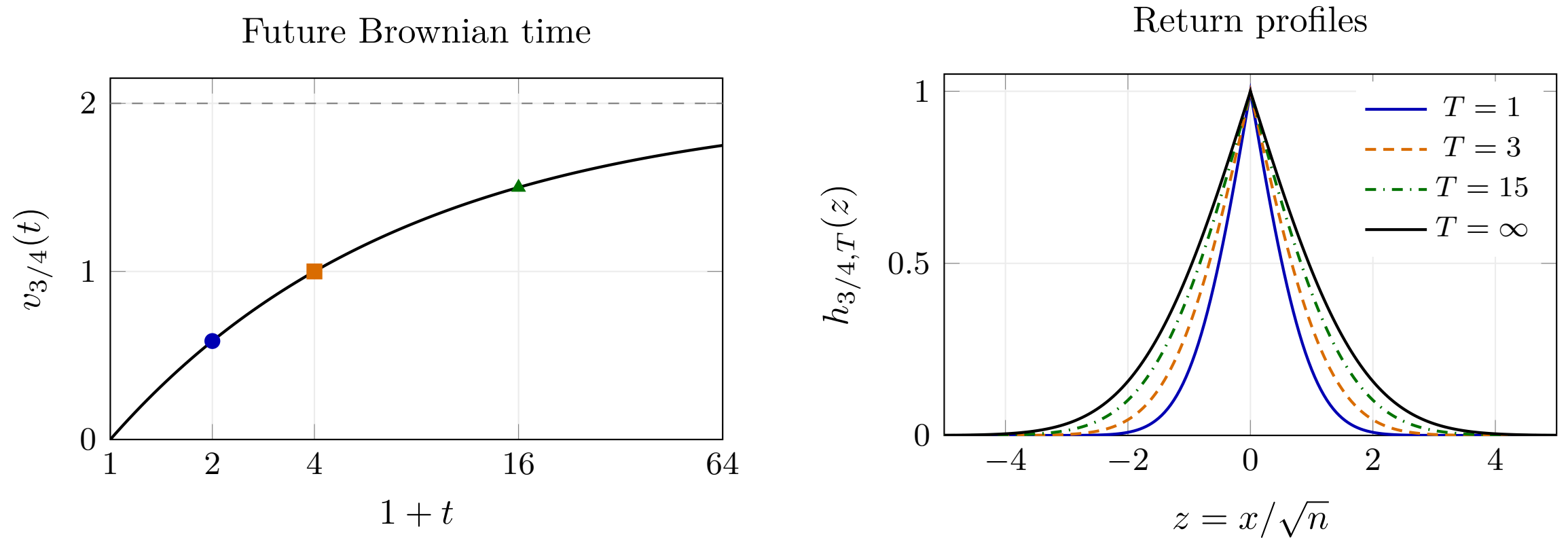
\begin{figure}[tbp]
\centering
\begin{minipage}{0.49\linewidth}
\centering
\begin{tikzpicture}
\begin{axis}[width=\linewidth,height=5.1cm,scale only axis=false,
 font=\small,xmode=log,log basis x=2,xmin=1,xmax=64,ymin=0,ymax=2.15,
 xtick={1,2,4,16,64},xticklabels={1,2,4,16,64},ytick={0,1,2},
 xlabel={$1+t$},ylabel={$v_{3/4}(t)$},grid=major,grid style={gray!15},
 title={Future Brownian time},title style={font=\small},
 tick label style={font=\footnotesize}]
\addplot[black,thick,domain=1:64,samples=100]{2*(1-1/sqrt(x))};
\addplot[gray,dashed,domain=1:64]{2};
\addplot[only marks,mark=*,blue!70!black] coordinates {(2,0.58578644)};
\addplot[only marks,mark=square*,orange!85!black] coordinates {(4,1)};
\addplot[only marks,mark=triangle*,green!45!black] coordinates {(16,1.5)};
\end{axis}
\end{tikzpicture}
\end{minipage}\hfill
\begin{minipage}{0.49\linewidth}
\centering
\begin{tikzpicture}
\begin{axis}[width=\linewidth,height=5.1cm,scale only axis=false,
 font=\small,xmin=-5,xmax=5,ymin=0,ymax=1.05,
 xtick={-4,-2,0,2,4},ytick={0,0.5,1},
 xlabel={$z=x/\sqrt n$},ylabel={$h_{3/4,T}(z)$},
 grid=major,grid style={gray!15},title={Return profiles},title style={font=\small},
 tick label style={font=\footnotesize},legend style={font=\scriptsize,
 at={(0.98,0.98)},anchor=north east,draw=none,fill=white,inner sep=2pt}]
\addplot[blue!70!black,solid,thick] coordinates {(-5.00,0.00000000) (-4.95,0.00000000) (-4.90,0.00000000) (-4.85,0.00000000) (-4.80,0.00000000) (-4.75,0.00000000) (-4.70,0.00000000) (-4.65,0.00000000) (-4.60,0.00000000) (-4.55,0.00000000) (-4.50,0.00000000) (-4.45,0.00000001) (-4.40,0.00000001) (-4.35,0.00000001) (-4.30,0.00000002) (-4.25,0.00000003) (-4.20,0.00000004) (-4.15,0.00000006) (-4.10,0.00000008) (-4.05,0.00000012) (-4.00,0.00000017) (-3.95,0.00000025) (-3.90,0.00000035) (-3.85,0.00000049) (-3.80,0.00000069) (-3.75,0.00000096) (-3.70,0.00000134) (-3.65,0.00000185) (-3.60,0.00000256) (-3.55,0.00000351) (-3.50,0.00000481) (-3.45,0.00000656) (-3.40,0.00000890) (-3.35,0.00001203) (-3.30,0.00001620) (-3.25,0.00002173) (-3.20,0.00002902) (-3.15,0.00003861) (-3.10,0.00005114) (-3.05,0.00006748) (-3.00,0.00008866) (-2.95,0.00011603) (-2.90,0.00015124) (-2.85,0.00019632) (-2.80,0.00025382) (-2.75,0.00032683) (-2.70,0.00041916) (-2.65,0.00053540) (-2.60,0.00068113) (-2.55,0.00086306) (-2.50,0.00108921) (-2.45,0.00136914) (-2.40,0.00171415) (-2.35,0.00213756) (-2.30,0.00265498) (-2.25,0.00328459) (-2.20,0.00404747) (-2.15,0.00496786) (-2.10,0.00607358) (-2.05,0.00739629) (-2.00,0.00897182) (-1.95,0.01084053) (-1.90,0.01304752) (-1.85,0.01564294) (-1.80,0.01868215) (-1.75,0.02222589) (-1.70,0.02634033) (-1.65,0.03109703) (-1.60,0.03657284) (-1.55,0.04284964) (-1.50,0.05001397) (-1.45,0.05815651) (-1.40,0.06737143) (-1.35,0.07775555) (-1.30,0.08940743) (-1.25,0.10242616) (-1.20,0.11691018) (-1.15,0.13295579) (-1.10,0.15065570) (-1.05,0.17009734) (-1.00,0.19136118) (-0.95,0.21451902) (-0.90,0.23963217) (-0.85,0.26674981) (-0.80,0.29590727) (-0.75,0.32712451) (-0.70,0.36040477) (-0.65,0.39573334) (-0.60,0.43307662) (-0.55,0.47238148) (-0.50,0.51357483) (-0.45,0.55656362) (-0.40,0.60123511) (-0.35,0.64745756) (-0.30,0.69508120) (-0.25,0.74393963) (-0.20,0.79385146) (-0.15,0.84462234) (-0.10,0.89604721) (-0.05,0.94791273) (0.00,1.00000000) (0.05,0.94791273) (0.10,0.89604721) (0.15,0.84462234) (0.20,0.79385146) (0.25,0.74393963) (0.30,0.69508120) (0.35,0.64745756) (0.40,0.60123511) (0.45,0.55656362) (0.50,0.51357483) (0.55,0.47238148) (0.60,0.43307662) (0.65,0.39573334) (0.70,0.36040477) (0.75,0.32712451) (0.80,0.29590727) (0.85,0.26674981) (0.90,0.23963217) (0.95,0.21451902) (1.00,0.19136118) (1.05,0.17009734) (1.10,0.15065570) (1.15,0.13295579) (1.20,0.11691018) (1.25,0.10242616) (1.30,0.08940743) (1.35,0.07775555) (1.40,0.06737143) (1.45,0.05815651) (1.50,0.05001397) (1.55,0.04284964) (1.60,0.03657284) (1.65,0.03109703) (1.70,0.02634033) (1.75,0.02222589) (1.80,0.01868215) (1.85,0.01564294) (1.90,0.01304752) (1.95,0.01084053) (2.00,0.00897182) (2.05,0.00739629) (2.10,0.00607358) (2.15,0.00496786) (2.20,0.00404747) (2.25,0.00328459) (2.30,0.00265498) (2.35,0.00213756) (2.40,0.00171415) (2.45,0.00136914) (2.50,0.00108921) (2.55,0.00086306) (2.60,0.00068113) (2.65,0.00053540) (2.70,0.00041916) (2.75,0.00032683) (2.80,0.00025382) (2.85,0.00019632) (2.90,0.00015124) (2.95,0.00011603) (3.00,0.00008866) (3.05,0.00006748) (3.10,0.00005114) (3.15,0.00003861) (3.20,0.00002902) (3.25,0.00002173) (3.30,0.00001620) (3.35,0.00001203) (3.40,0.00000890) (3.45,0.00000656) (3.50,0.00000481) (3.55,0.00000351) (3.60,0.00000256) (3.65,0.00000185) (3.70,0.00000134) (3.75,0.00000096) (3.80,0.00000069) (3.85,0.00000049) (3.90,0.00000035) (3.95,0.00000025) (4.00,0.00000017) (4.05,0.00000012) (4.10,0.00000008) (4.15,0.00000006) (4.20,0.00000004) (4.25,0.00000003) (4.30,0.00000002) (4.35,0.00000001) (4.40,0.00000001) (4.45,0.00000001) (4.50,0.00000000) (4.55,0.00000000) (4.60,0.00000000) (4.65,0.00000000) (4.70,0.00000000) (4.75,0.00000000) (4.80,0.00000000) (4.85,0.00000000) (4.90,0.00000000) (4.95,0.00000000) (5.00,0.00000000)};
\addlegendentry{$T=1$}
\addplot[orange!85!black,densely dashed,thick] coordinates {(-5.00,0.00000057) (-4.95,0.00000074) (-4.90,0.00000096) (-4.85,0.00000123) (-4.80,0.00000159) (-4.75,0.00000203) (-4.70,0.00000260) (-4.65,0.00000332) (-4.60,0.00000422) (-4.55,0.00000536) (-4.50,0.00000680) (-4.45,0.00000859) (-4.40,0.00001083) (-4.35,0.00001361) (-4.30,0.00001708) (-4.25,0.00002138) (-4.20,0.00002669) (-4.15,0.00003325) (-4.10,0.00004132) (-4.05,0.00005122) (-4.00,0.00006334) (-3.95,0.00007815) (-3.90,0.00009619) (-3.85,0.00011812) (-3.80,0.00014470) (-3.75,0.00017683) (-3.70,0.00021560) (-3.65,0.00026224) (-3.60,0.00031822) (-3.55,0.00038523) (-3.50,0.00046526) (-3.45,0.00056059) (-3.40,0.00067386) (-3.35,0.00080812) (-3.30,0.00096685) (-3.25,0.00115405) (-3.20,0.00137428) (-3.15,0.00163270) (-3.10,0.00193521) (-3.05,0.00228841) (-3.00,0.00269980) (-2.95,0.00317774) (-2.90,0.00373163) (-2.85,0.00437192) (-2.80,0.00511026) (-2.75,0.00595953) (-2.70,0.00693395) (-2.65,0.00804918) (-2.60,0.00932238) (-2.55,0.01077229) (-2.50,0.01241933) (-2.45,0.01428562) (-2.40,0.01639507) (-2.35,0.01877341) (-2.30,0.02144822) (-2.25,0.02444895) (-2.20,0.02780690) (-2.15,0.03155521) (-2.10,0.03572884) (-2.05,0.04036443) (-2.00,0.04550026) (-1.95,0.05117612) (-1.90,0.05743312) (-1.85,0.06431355) (-1.80,0.07186064) (-1.75,0.08011831) (-1.70,0.08913093) (-1.65,0.09894294) (-1.60,0.10959858) (-1.55,0.12114152) (-1.50,0.13361440) (-1.45,0.14705852) (-1.40,0.16151332) (-1.35,0.17701598) (-1.30,0.19360097) (-1.25,0.21129955) (-1.20,0.23013934) (-1.15,0.25014387) (-1.10,0.27133212) (-1.05,0.29371811) (-1.00,0.31731051) (-0.95,0.34211225) (-0.90,0.36812025) (-0.85,0.39532509) (-0.80,0.42371080) (-0.75,0.45325470) (-0.70,0.48392730) (-0.65,0.51569222) (-0.60,0.54850624) (-0.55,0.58231937) (-0.50,0.61707508) (-0.45,0.65271044) (-0.40,0.68915652) (-0.35,0.72633870) (-0.30,0.76417716) (-0.25,0.80258735) (-0.20,0.84148058) (-0.15,0.88076462) (-0.10,0.92034433) (-0.05,0.96012239) (0.00,1.00000000) (0.05,0.96012239) (0.10,0.92034433) (0.15,0.88076462) (0.20,0.84148058) (0.25,0.80258735) (0.30,0.76417716) (0.35,0.72633870) (0.40,0.68915652) (0.45,0.65271044) (0.50,0.61707508) (0.55,0.58231937) (0.60,0.54850624) (0.65,0.51569222) (0.70,0.48392730) (0.75,0.45325470) (0.80,0.42371080) (0.85,0.39532509) (0.90,0.36812025) (0.95,0.34211225) (1.00,0.31731051) (1.05,0.29371811) (1.10,0.27133212) (1.15,0.25014387) (1.20,0.23013934) (1.25,0.21129955) (1.30,0.19360097) (1.35,0.17701598) (1.40,0.16151332) (1.45,0.14705852) (1.50,0.13361440) (1.55,0.12114152) (1.60,0.10959858) (1.65,0.09894294) (1.70,0.08913093) (1.75,0.08011831) (1.80,0.07186064) (1.85,0.06431355) (1.90,0.05743312) (1.95,0.05117612) (2.00,0.04550026) (2.05,0.04036443) (2.10,0.03572884) (2.15,0.03155521) (2.20,0.02780690) (2.25,0.02444895) (2.30,0.02144822) (2.35,0.01877341) (2.40,0.01639507) (2.45,0.01428562) (2.50,0.01241933) (2.55,0.01077229) (2.60,0.00932238) (2.65,0.00804918) (2.70,0.00693395) (2.75,0.00595953) (2.80,0.00511026) (2.85,0.00437192) (2.90,0.00373163) (2.95,0.00317774) (3.00,0.00269980) (3.05,0.00228841) (3.10,0.00193521) (3.15,0.00163270) (3.20,0.00137428) (3.25,0.00115405) (3.30,0.00096685) (3.35,0.00080812) (3.40,0.00067386) (3.45,0.00056059) (3.50,0.00046526) (3.55,0.00038523) (3.60,0.00031822) (3.65,0.00026224) (3.70,0.00021560) (3.75,0.00017683) (3.80,0.00014470) (3.85,0.00011812) (3.90,0.00009619) (3.95,0.00007815) (4.00,0.00006334) (4.05,0.00005122) (4.10,0.00004132) (4.15,0.00003325) (4.20,0.00002669) (4.25,0.00002138) (4.30,0.00001708) (4.35,0.00001361) (4.40,0.00001083) (4.45,0.00000859) (4.50,0.00000680) (4.55,0.00000536) (4.60,0.00000422) (4.65,0.00000332) (4.70,0.00000260) (4.75,0.00000203) (4.80,0.00000159) (4.85,0.00000123) (4.90,0.00000096) (4.95,0.00000074) (5.00,0.00000057)};
\addlegendentry{$T=3$}
\addplot[green!45!black,dashdotted,thick] coordinates {(-5.00,0.00004456) (-4.95,0.00005307) (-4.90,0.00006312) (-4.85,0.00007495) (-4.80,0.00008885) (-4.75,0.00010516) (-4.70,0.00012428) (-4.65,0.00014663) (-4.60,0.00017273) (-4.55,0.00020316) (-4.50,0.00023856) (-4.45,0.00027970) (-4.40,0.00032741) (-4.35,0.00038266) (-4.30,0.00044653) (-4.25,0.00052024) (-4.20,0.00060517) (-4.15,0.00070286) (-4.10,0.00081504) (-4.05,0.00094364) (-4.00,0.00109084) (-3.95,0.00125902) (-3.90,0.00145086) (-3.85,0.00166934) (-3.80,0.00191772) (-3.75,0.00219965) (-3.70,0.00251910) (-3.65,0.00288048) (-3.60,0.00328861) (-3.55,0.00374877) (-3.50,0.00426672) (-3.45,0.00484876) (-3.40,0.00550173) (-3.35,0.00623303) (-3.30,0.00705073) (-3.25,0.00796349) (-3.20,0.00898067) (-3.15,0.01011233) (-3.10,0.01136927) (-3.05,0.01276301) (-3.00,0.01430588) (-2.95,0.01601099) (-2.90,0.01789227) (-2.85,0.01996445) (-2.80,0.02224312) (-2.75,0.02474467) (-2.70,0.02748634) (-2.65,0.03048616) (-2.60,0.03376298) (-2.55,0.03733642) (-2.50,0.04122683) (-2.45,0.04545529) (-2.40,0.05004352) (-2.35,0.05501383) (-2.30,0.06038908) (-2.25,0.06619258) (-2.20,0.07244801) (-2.15,0.07917935) (-2.10,0.08641073) (-2.05,0.09416638) (-2.00,0.10247043) (-1.95,0.11134689) (-1.90,0.12081939) (-1.85,0.13091113) (-1.80,0.14164469) (-1.75,0.15304188) (-1.70,0.16512359) (-1.65,0.17790959) (-1.60,0.19141843) (-1.55,0.20566718) (-1.50,0.22067136) (-1.45,0.23644472) (-1.40,0.25299906) (-1.35,0.27034414) (-1.30,0.28848746) (-1.25,0.30743417) (-1.20,0.32718688) (-1.15,0.34774560) (-1.10,0.36910759) (-1.05,0.39126728) (-1.00,0.41421618) (-0.95,0.43794282) (-0.90,0.46243273) (-0.85,0.48766835) (-0.80,0.51362911) (-0.75,0.54029137) (-0.70,0.56762850) (-0.65,0.59561090) (-0.60,0.62420611) (-0.55,0.65337891) (-0.50,0.68309140) (-0.45,0.71330317) (-0.40,0.74397148) (-0.35,0.77505137) (-0.30,0.80649594) (-0.25,0.83825649) (-0.20,0.87028277) (-0.15,0.90252325) (-0.10,0.93492531) (-0.05,0.96743555) (0.00,1.00000000) (0.05,0.96743555) (0.10,0.93492531) (0.15,0.90252325) (0.20,0.87028277) (0.25,0.83825649) (0.30,0.80649594) (0.35,0.77505137) (0.40,0.74397148) (0.45,0.71330317) (0.50,0.68309140) (0.55,0.65337891) (0.60,0.62420611) (0.65,0.59561090) (0.70,0.56762850) (0.75,0.54029137) (0.80,0.51362911) (0.85,0.48766835) (0.90,0.46243273) (0.95,0.43794282) (1.00,0.41421618) (1.05,0.39126728) (1.10,0.36910759) (1.15,0.34774560) (1.20,0.32718688) (1.25,0.30743417) (1.30,0.28848746) (1.35,0.27034414) (1.40,0.25299906) (1.45,0.23644472) (1.50,0.22067136) (1.55,0.20566718) (1.60,0.19141843) (1.65,0.17790959) (1.70,0.16512359) (1.75,0.15304188) (1.80,0.14164469) (1.85,0.13091113) (1.90,0.12081939) (1.95,0.11134689) (2.00,0.10247043) (2.05,0.09416638) (2.10,0.08641073) (2.15,0.07917935) (2.20,0.07244801) (2.25,0.06619258) (2.30,0.06038908) (2.35,0.05501383) (2.40,0.05004352) (2.45,0.04545529) (2.50,0.04122683) (2.55,0.03733642) (2.60,0.03376298) (2.65,0.03048616) (2.70,0.02748634) (2.75,0.02474467) (2.80,0.02224312) (2.85,0.01996445) (2.90,0.01789227) (2.95,0.01601099) (3.00,0.01430588) (3.05,0.01276301) (3.10,0.01136927) (3.15,0.01011233) (3.20,0.00898067) (3.25,0.00796349) (3.30,0.00705073) (3.35,0.00623303) (3.40,0.00550173) (3.45,0.00484876) (3.50,0.00426672) (3.55,0.00374877) (3.60,0.00328861) (3.65,0.00288048) (3.70,0.00251910) (3.75,0.00219965) (3.80,0.00191772) (3.85,0.00166934) (3.90,0.00145086) (3.95,0.00125902) (4.00,0.00109084) (4.05,0.00094364) (4.10,0.00081504) (4.15,0.00070286) (4.20,0.00060517) (4.25,0.00052024) (4.30,0.00044653) (4.35,0.00038266) (4.40,0.00032741) (4.45,0.00027970) (4.50,0.00023856) (4.55,0.00020316) (4.60,0.00017273) (4.65,0.00014663) (4.70,0.00012428) (4.75,0.00010516) (4.80,0.00008885) (4.85,0.00007495) (4.90,0.00006312) (4.95,0.00005307) (5.00,0.00004456)};
\addlegendentry{$T=15$}
\addplot[black,solid,thick] coordinates {(-5.00,0.00040695) (-4.95,0.00046495) (-4.90,0.00053058) (-4.85,0.00060477) (-4.80,0.00068851) (-4.75,0.00078294) (-4.70,0.00088927) (-4.65,0.00100885) (-4.60,0.00114318) (-4.55,0.00129387) (-4.50,0.00146272) (-4.45,0.00165167) (-4.40,0.00186285) (-4.35,0.00209858) (-4.30,0.00236139) (-4.25,0.00265403) (-4.20,0.00297947) (-4.15,0.00334093) (-4.10,0.00374190) (-4.05,0.00418615) (-4.00,0.00467773) (-3.95,0.00522101) (-3.90,0.00582067) (-3.85,0.00648173) (-3.80,0.00720957) (-3.75,0.00800994) (-3.70,0.00888897) (-3.65,0.00985318) (-3.60,0.01090950) (-3.55,0.01206529) (-3.50,0.01332833) (-3.45,0.01470685) (-3.40,0.01620954) (-3.35,0.01784553) (-3.30,0.01962441) (-3.25,0.02155627) (-3.20,0.02365162) (-3.15,0.02592146) (-3.10,0.02837727) (-3.05,0.03103094) (-3.00,0.03389485) (-2.95,0.03698179) (-2.90,0.04030497) (-2.85,0.04387801) (-2.80,0.04771488) (-2.75,0.05182993) (-2.70,0.05623780) (-2.65,0.06095345) (-2.60,0.06599206) (-2.55,0.07136901) (-2.50,0.07709987) (-2.45,0.08320030) (-2.40,0.08968602) (-2.35,0.09657275) (-2.30,0.10387616) (-2.25,0.11161177) (-2.20,0.11979493) (-2.15,0.12844072) (-2.10,0.13756389) (-2.05,0.14717877) (-2.00,0.15729921) (-1.95,0.16793847) (-1.90,0.17910919) (-1.85,0.19082327) (-1.80,0.20309179) (-1.75,0.21592494) (-1.70,0.22933194) (-1.65,0.24332096) (-1.60,0.25789904) (-1.55,0.27307199) (-1.50,0.28884437) (-1.45,0.30521938) (-1.40,0.32219881) (-1.35,0.33978297) (-1.30,0.35797067) (-1.25,0.37675912) (-1.20,0.39614391) (-1.15,0.41611899) (-1.10,0.43667663) (-1.05,0.45780739) (-1.00,0.47950012) (-0.95,0.50174195) (-0.90,0.52451828) (-0.85,0.54781284) (-0.80,0.57160764) (-0.75,0.59588309) (-0.70,0.62061795) (-0.65,0.64578943) (-0.60,0.67137324) (-0.55,0.69734366) (-0.50,0.72367361) (-0.45,0.75033471) (-0.40,0.77729741) (-0.35,0.80453106) (-0.30,0.83200403) (-0.25,0.85968380) (-0.20,0.88753708) (-0.15,0.91552997) (-0.10,0.94362802) (-0.05,0.97179640) (0.00,1.00000000) (0.05,0.97179640) (0.10,0.94362802) (0.15,0.91552997) (0.20,0.88753708) (0.25,0.85968380) (0.30,0.83200403) (0.35,0.80453106) (0.40,0.77729741) (0.45,0.75033471) (0.50,0.72367361) (0.55,0.69734366) (0.60,0.67137324) (0.65,0.64578943) (0.70,0.62061795) (0.75,0.59588309) (0.80,0.57160764) (0.85,0.54781284) (0.90,0.52451828) (0.95,0.50174195) (1.00,0.47950012) (1.05,0.45780739) (1.10,0.43667663) (1.15,0.41611899) (1.20,0.39614391) (1.25,0.37675912) (1.30,0.35797067) (1.35,0.33978297) (1.40,0.32219881) (1.45,0.30521938) (1.50,0.28884437) (1.55,0.27307199) (1.60,0.25789904) (1.65,0.24332096) (1.70,0.22933194) (1.75,0.21592494) (1.80,0.20309179) (1.85,0.19082327) (1.90,0.17910919) (1.95,0.16793847) (2.00,0.15729921) (2.05,0.14717877) (2.10,0.13756389) (2.15,0.12844072) (2.20,0.11979493) (2.25,0.11161177) (2.30,0.10387616) (2.35,0.09657275) (2.40,0.08968602) (2.45,0.08320030) (2.50,0.07709987) (2.55,0.07136901) (2.60,0.06599206) (2.65,0.06095345) (2.70,0.05623780) (2.75,0.05182993) (2.80,0.04771488) (2.85,0.04387801) (2.90,0.04030497) (2.95,0.03698179) (3.00,0.03389485) (3.05,0.03103094) (3.10,0.02837727) (3.15,0.02592146) (3.20,0.02365162) (3.25,0.02155627) (3.30,0.01962441) (3.35,0.01784553) (3.40,0.01620954) (3.45,0.01470685) (3.50,0.01332833) (3.55,0.01206529) (3.60,0.01090950) (3.65,0.00985318) (3.70,0.00888897) (3.75,0.00800994) (3.80,0.00720957) (3.85,0.00648173) (3.90,0.00582067) (3.95,0.00522101) (4.00,0.00467773) (4.05,0.00418615) (4.10,0.00374190) (4.15,0.00334093) (4.20,0.00297947) (4.25,0.00265403) (4.30,0.00236139) (4.35,0.00209858) (4.40,0.00186285) (4.45,0.00165167) (4.50,0.00146272) (4.55,0.00129387) (4.60,0.00114318) (4.65,0.00100885) (4.70,0.00088927) (4.75,0.00078294) (4.80,0.00068851) (4.85,0.00060477) (4.90,0.00053058) (4.95,0.00046495) (5.00,0.00040695)};
\addlegendentry{$T=\infty$}
\end{axis}
\end{tikzpicture}
\end{minipage}
\caption{Future variance and return profiles for $p=7/8$, hence $a=3/4$.
Left: $v_a(t)$ approaches $2$ as the physical horizon $n(1+t)$ tends to infinity; the time axis is logarithmic. The marked horizons $T=1,3,15$ correspond to the curves on the right, which approach the eventual-return profile $h_a$. Its full integral over $\R$ is $\kappa_a=4/\sqrt\pi$, the factor multiplying $f_a(0)$ in the last-return asymptotic. The curves show the limiting formulas, not finite-$n$ simulations.}
\label{fig:clock-kernel}
\end{figure}

The deterministic clock used below is
\begin{equation}\label{eq:clock}
 v_a(t)=\int_0^t(1+s)^{-2a}\dd s
       =\frac{1-(1+t)^{-\delta}}\delta,\qquad v_a(\infty)=\frac1\delta.
\end{equation}
The clock $v_a(t)$ is the limiting variance accumulated over the elapsed time $nt$, after division by $n$. Later increments receive smaller weights, whose squares behave like $(1+t)^{-2a}$. Since $2a>1$, infinitely many future steps contribute only the finite Brownian time $1/\delta$. This is why the return kernel below involves Brownian motion stopped at a fixed horizon; Figure~\ref{fig:clock-kernel} illustrates the clock and its hitting profiles.

\citet[Theorem 1.1]{Tok} obtains the same Brownian clock under the original law. We need a limit under $\Pp_{n,x_n}$, with deterministic $x_n=O(\sqrt n)$. Convergence under the original law does not directly imply this limit, since conditioning on reachable states involves probabilities tending to zero. We verify the hypotheses of \citet[Theorem 2.1(ii)]{Whitt}, using Lemma \ref{lem:second} to control the random correction to the quadratic variation.

\begin{theorem}\label{thm:future-fclt}
Let $x_n$ be admissible at time $n$ and $x_n/\sqrt n\to z\in\R$. Under $\Pp_{n,x_n}$,
\begin{equation}\label{eq:future-fclt}
 \left(\frac{R_{n,n+\lfloor nt\rfloor}S_{n+\lfloor nt\rfloor}-x_n}{\sqrt n}\right)_{t\geq0}
 \law (B(v_a(t)))_{t\geq0}
\end{equation}
in $D([0,\infty))$ with the local $J_1$ topology, where $B$ is a standard Brownian motion.
\end{theorem}
\begin{proof}
Fix $0<T<\infty$ and denote the process on the left by $Y_n(t)$. Since $R_{n,m}\leq1$, \eqref{eq:jumps} bounds its jumps by $2/\sqrt n$. Formula \eqref{eq:incrementvar} gives the quadratic-variation decomposition from \citet[(A.7)]{Bercu}, with the present starting time and normalization:
\begin{align*}
 \langle Y_n\rangle_t&=\frac1n\sum_{m=n}^{n+\lfloor nt\rfloor-1}R_{n,m+1}^2-D_n(t),\\
 D_n(t)&=\frac{a^2}n\sum_{m=n}^{n+\lfloor nt\rfloor-1}
          R_{n,m+1}^2\frac{S_m^2}{m^2}.
\end{align*}
By \eqref{eq:Rlimit}, the first term converges uniformly on $[0,T]$ to $v_a(t)$ as a Riemann sum. The second is nonnegative and increasing. Lemma \ref{lem:second} gives
\[
 \E_{n,x_n}[D_n(T)]\leq\frac Cn
 \sum_{m=n}^{n+\lfloor nT\rfloor-1}\frac n{m^2}
 \leq\frac{C(T+1)}n\longrightarrow0.
\]
Hence $\langle Y_n\rangle\inprob v_a$ uniformly on $[0,T]$. To apply the finite-interval form of the martingale functional central limit theorem in \citet[Theorem 2.1(ii)]{Whitt}, put $\theta=v_a^{-1}$ and $\widetilde Y_n(u)=Y_n(\theta(u))$, $0\leq u\leq v_a(T)$. Then
\[
 \sup_{u\leq v_a(T)}|\langle\widetilde Y_n\rangle_u-u|\inprob0,
 \qquad
 \sup_u|\Delta\widetilde Y_n(u)|\leq\frac2{\sqrt n},
 \qquad
 \sup_u\Delta\langle\widetilde Y_n\rangle_u\leq\frac4n.
\]
The last inequality follows because jumps occur at deterministic times and bracket jumps are conditional second moments. Thus the bracket and small-jump hypotheses of the cited theorem hold, and $\widetilde Y_n\law B$. Composing with $v_a$ gives \eqref{eq:future-fclt} on $[0,T]$, and $T$ is arbitrary.
\end{proof}

\section{Probability of a future return}\label{sec:hitting}
For admissible $(n,x)$, put
\[
 \tau_0^{(n)}=\inf\{m\geq n:S_m=0\},\quad
 H_n(x)=\Pp_{n,x}(\tau_0^{(n)}<\infty),\quad
 H_{n,T}(x)=\Pp_{n,x}(\tau_0^{(n)}\leq n+\lfloor nT\rfloor).
\]
In particular, a walk started at zero has hitting time $n$ in these definitions.

Define
\begin{equation}\label{eq:ha}
 h_a(z)=2\Phi(-\sqrt\delta\,|z|),\qquad
 \kappa_a=\int_\R h_a(z)\dd z=\sqrt{\frac8{\pi\delta}},
\end{equation}
where $\Phi$ is the standard normal distribution function. By the reflection principle, $h_a(z)$ is the probability that Brownian motion started at $z$ hits zero by time $1/\delta$; see \citet[Theorems~2.19--2.21]{MP}. It weights the scaled starting positions according to their chance of returning. Its area $\kappa_a$ converts the nearly constant local mass supplied by the LLT into the last-return constant $f_a(0)\kappa_a$. The following theorem makes this approximation uniform and controls its spatial tails.

\begin{theorem}\label{thm:kernel}
For every $0<K<\infty$,
\[
 \sup_{\substack{|x|\leq K\sqrt n\\(n,x)\text{ admissible}}}
 \left|H_n(x)-h_a(x/\sqrt n)\right|\longrightarrow0.
\]
There are constants $c_a,C_a>0$ such that, for every admissible $(n,x)$,
\begin{equation}\label{eq:gaussian-intro}
 H_n(x)\leq C_a\exp(-c_a x^2/n).
\end{equation}
\end{theorem}

The proof has three steps. Functional convergence gives the return probability on a fixed time interval. A late-hit estimate then permits the horizon to tend to infinity, and an exponential bound controls distant starting points. These last two bounds address time and space separately. For the first step we use the minimum of the normalized path; compare \citet[proof of Theorem 1.4]{FangTraining}.

Define also
\begin{equation}\label{eq:finitekernel}
 h_{a,T}(z)=2\Phi\left(-\frac{|z|}{\sqrt{v_a(T)}}\right),\qquad T>0.
\end{equation}

\begin{lemma}\label{lem:finite-hitting}
If $x_n$ is admissible at time $n$ and $x_n/\sqrt n\to z\in\R$, then $H_{n,T}(x_n)\to h_{a,T}(z)$ for every $T>0$.
\end{lemma}
\begin{proof}
Reflection of the transition rule gives $H_{n,T}(x)=H_{n,T}(-x)$, so it suffices to consider $x_n\geq0$ and $z\geq0$. By Theorem \ref{thm:future-fclt},
\[
 Z_n(t):=\frac{R_{n,n+\lfloor nt\rfloor}S_{n+\lfloor nt\rfloor}}{\sqrt n}
 \law z+B(v_a(t)),\qquad 0\leq t\leq T.
\]
Since $R_{n,m}>0$ and the walk cannot cross zero without visiting it,
\[
 \{\tau_0^{(n)}\leq n+\lfloor nT\rfloor\}
 =\left\{\inf_{0\leq t\leq T}Z_n(t)\leq0\right\}.
\]
This includes $x_n=0$, since the initial visit is counted. The minimum is continuous in the $J_1$ topology at continuous paths, because convergence to such a path is uniform on $[0,T]$. By the Brownian reflection principle \citep[Chapter 2]{MP}, the limiting minimum has the law of $z-|G_T|$, where $G_T\sim N(0,v_a(T))$; it has no atom at zero, including when $z=0$. The continuous-mapping and continuity-set implications of the Portmanteau theorem \citep[Theorem 3.10.1]{Durrett} therefore give
\[
 H_{n,T}(x_n)\longrightarrow
 \Pp\left(\inf_{s\leq v_a(T)}(z+B_s)\leq0\right)
 =2\Phi\left(-\frac z{\sqrt{v_a(T)}}\right).
\]
\end{proof}

To pass from finite to infinite horizons, we need to control first hits after a large rescaled time, even though the Brownian clock has the finite endpoint $v_a(\infty)=1/\delta$.

\begin{lemma}\label{lem:late-hitting}
For every $0<K<\infty$,
\[
 \lim_{T\to\infty}\limsup_{n\to\infty}
 \sup_{\substack{|x|\leq K\sqrt n\\(n,x)\text{ admissible}}}
 \Pp_{n,x}(n+\lfloor nT\rfloor<\tau_0^{(n)}<\infty)=0.
\]
\end{lemma}
\begin{proof}
Let $m_T=n+\lfloor nT\rfloor$ and $\eta>0$. A hit after $m_T$ implies either $|M_{m_T}^{(n)}|\leq\eta\sqrt n$ or $\sup_{m\geq m_T}|M_m^{(n)}-M_{m_T}^{(n)}|\geq\eta\sqrt n$. Orthogonality of martingale increments and \eqref{eq:jumps} give, for $r>m_T$,
\[
 \E_{n,x}[(M_r^{(n)}-M_{m_T}^{(n)})^2]
 \leq4\sum_{m\geq m_T}R_{n,m+1}^2<\infty.
\]
The martingale converges in $L^2$, so Doob's $L^2$ maximal inequality applies up to the infinite horizon; see \citet[Theorems 4.4.6 and 4.4.4]{Durrett}. Together with \eqref{eq:variance-tail}, it gives
\begin{equation}\label{eq:Doob-tail}
 \Pp_{n,x}\left(\sup_{m\geq m_T}|M_m^{(n)}-M_{m_T}^{(n)}|\geq\eta\sqrt n\right)
 \leq\frac{16}{\eta^2n}\sum_{m\geq m_T}R_{n,m+1}^2
 \leq C_a\eta^{-2}(1+T)^{-\delta}.
\end{equation}
For fixed $T\geq1$, Theorem \ref{thm:future-fclt} implies
\begin{equation}\label{eq:small-interval}
 \limsup_n\sup_{|x|\leq K\sqrt n}
 \Pp_{n,x}(|M_{m_T}^{(n)}|\leq\eta\sqrt n)\leq C_a\eta,
\end{equation}
where only admissible states are included. Indeed, any sequence of starting points in this window has a subsequence with $x_n/\sqrt n\to z$. Along it, the limit is $N(z,v_a(T))$, whose density is bounded by $\{2\pi v_a(1)\}^{-1/2}$, uniformly in $z$ and $T\geq1$. Applying the closed-set bound in the Portmanteau theorem \citep[Theorem 3.10.1(iii)]{Durrett} to $[-\eta,\eta]$ proves \eqref{eq:small-interval}.

Combining \eqref{eq:Doob-tail} and \eqref{eq:small-interval} bounds the limsup in the lemma by
\[
 C_a\{\eta+\eta^{-2}(1+T)^{-\delta}\}.
\]
Choose $\eta=(1+T)^{-\delta/4}$ and let $T\to\infty$.
\end{proof}

\begin{proof}[Proof of the convergence assertion in Theorem \ref{thm:kernel}]
Let $x_n$ be admissible at time $n$ and $x_n/\sqrt n\to z\in\R$. For every $T>0$, Lemma \ref{lem:finite-hitting} gives $H_{n,T}(x_n)\to h_{a,T}(z)$, while
\[
 0\leq H_n(x_n)-H_{n,T}(x_n)
 =\Pp_{n,x_n}(n+\lfloor nT\rfloor<\tau_0^{(n)}<\infty).
\]
The scaled starting points remain bounded, so Lemma \ref{lem:late-hitting} makes the limsup of this difference tend to zero as $T\to\infty$. Since $h_{a,T}(z)\to h_a(z)$, we obtain $H_n(x_n)\to h_a(z)$.
Uniformity on a bounded scaled window follows by compactness: a sequence of errors bounded away from zero would have a subsequence with $x_n/\sqrt n\to z$, contradicting this sequential convergence and continuity of $h_a$.
\end{proof}

The global spatial bound uses the Azuma--Hoeffding exponential-martingale method \citep{Hoeffding,Azuma}. \citet[Lemma 3.1]{FangTraining} also uses this method to localize a first return for trained diffusive ERWs. Here the summability of squared future increments gives a bound over the entire future.

\begin{lemma}\label{lem:gaussian}
The bound \eqref{eq:gaussian-intro} holds.
\end{lemma}
\begin{proof}
For a martingale $N_j$ with deterministic bounds $|\Delta N_j|\leq c_j$, we use the maximal Azuma--Hoeffding bound \citep{Hoeffding,Azuma}, recalling its exponential-supermartingale proof:
\begin{equation}\label{eq:Azuma}
 \Pp\left(\max_{j\leq r}|N_j-N_0|\geq u\right)
 \leq2\exp\left(-\frac{u^2}{2\sum_{j=1}^rc_j^2}\right).
\end{equation}
Indeed, the conditional form of Hoeffding's exponential-moment bound \citep[(4.16)]{Hoeffding} makes $\exp\{\lambda(N_j-N_0)-\lambda^2\sum_{k\leq j}c_k^2/2\}$ a supermartingale. Stopping it at its first crossing of a fixed level gives the maximal bound; optimizing at $\lambda=u/\sum_{j\leq r}c_j^2$ for each sign yields \eqref{eq:Azuma}.
For $M_m^{(n)}$, Lemmas \ref{lem:normalization} and \ref{lem:martingale} give
\[
 \sum_{m=n}^\infty(2R_{n,m+1})^2
 \leq C_a n^{2a}\sum_{m=n}^\infty(m+1)^{-2a}\leq C_an.
\]
A hit of zero requires a displacement of size $|x|$ from $M_n^{(n)}=x$. Apply \eqref{eq:Azuma} to a finite horizon and then let that horizon increase to infinity. This proves $H_n(x)\leq2\exp(-x^2/(C_an))$. The case $x=0$ is included by enlarging the prefactor.
\end{proof}

\begin{proposition}\label{prop:finite-uniform}
For $T>0$ and $0<K<\infty$,
\[
 \sup_{\substack{|x|\leq K\sqrt n\\(n,x)\text{ admissible}}}
 |H_{n,T}(x)-h_{a,T}(x/\sqrt n)|\longrightarrow0.
\]
\end{proposition}
\begin{proof}
Lemma \ref{lem:finite-hitting} gives convergence along every sequence of scaled starting points. Compactness of $[-K,K]$ and continuity of $h_{a,T}$ give uniformity by the same argument as above.
\end{proof}

The symmetry of the transition also gives an exact relation between hitting and one-time sign probabilities.

\begin{lemma}\label{lem:reflection}
For every admissible $(n,x)$, $M_m^{(n)}$ has an almost-sure and $L^2$ limit $M_\infty^{(n)}$, with $\Pp_{n,x}(M_\infty^{(n)}=0)=0$. If $x\ne0$ and $m\geq n$, then
\begin{align}
 \Pp_{n,x}(\tau_0^{(n)}\leq m)
 &=2\Pp_{n,x}(xS_m<0)+\Pp_{n,x}(S_m=0),\label{eq:discrete-reflection}\\
 H_n(x)&=2\Pp_{n,x}(xM_\infty^{(n)}<0).\label{eq:terminal-reflection}
\end{align}
\end{lemma}
\begin{proof}
By Lemmas \ref{lem:normalization} and \ref{lem:martingale}, the sum of the second moments of the martingale increments is finite. Hence the martingale converges almost surely and in $L^2$. Under the original law with $q=1/2$, every admissible event $\{S_n=x\}$ has positive probability, since $0<p<1$. Conditioning the almost-sure convergence and nonvanishing limit in Proposition \ref{thm:llt} on this event transfers these properties to $\Pp_{n,x}$. The product formula gives
\[
 M_\infty^{(n)}=\frac{\Gamma(n+a)}{\Gamma(n)}\lim_{k\to\infty}\frac{S_k}{k^a}\ne0
 \quad\Pp_{n,x}\text{-almost surely}.
\]
In particular its sign is the eventual sign of $S_k$.

Suppose $x>0$. A negative position at time $m$ requires a previous zero. Conditional on $\tau_0^{(n)}=k\leq m$, the future from $(k,0)$ is invariant under reflection, because \eqref{eq:transition} is sign-symmetric. Thus the positive and negative endpoint probabilities after this hit are equal. Summing over $k$ and including the zero endpoint gives \eqref{eq:discrete-reflection}. The same reflection exchanges the eventual signs, giving \eqref{eq:terminal-reflection}. The case $x<0$ follows by reflection.
\end{proof}

These identities also give an alternative derivation of the limiting kernels from scalar limits. For $m_T=n+\lfloor nT\rfloor$, \eqref{eq:incrementvar} and \eqref{eq:variance-tail} imply
\[
 \frac1n\E_{n,x}\bigl[(M_\infty^{(n)}-M_{m_T}^{(n)})^2\bigr]
 \leq C_a(1+T)^{-\delta}.
\]
Together with Theorem \ref{thm:future-fclt}, this yields $M_\infty^{(n)}/\sqrt n\law N(z,1/\delta)$ under $\Pp_{n,x_n}$ whenever $x_n/\sqrt n\to z$, by taking $n\to\infty$ and then $T\to\infty$. Equation \eqref{eq:terminal-reflection} recovers $h_a(z)$; \eqref{eq:discrete-reflection} and the fixed-time limit similarly recover $h_{a,T}(z)$, since the limiting normal has no atom at zero. At $z=0$, split into subsequences with $x_n>0$, $x_n<0$ or $x_n=0$; the latter have hitting probability one.

\section{The tail of the last return}\label{sec:tail}
The tail exponent is $a-1/2=\delta/2$. We combine the LLT of \citet[Theorem 3.1]{TV}, recorded in Proposition \ref{thm:llt}, with the hitting estimates of Theorem \ref{thm:kernel}, first on bounded $\sqrt n$ windows and then over the full lattice.
Conditioning on $S_n$ gives the exact decomposition
\begin{equation}\label{eq:decomposition}
 \Pp(E_n)=\sum_{\substack{x\in\Z\\x\equiv n\ (\mathrm{mod}\ 2)}}
 \Pp(S_n=x)H_n(x).
\end{equation}

\begin{lemma}\label{lem:central}
For every $0<K<\infty$,
\[
 n^{\delta/2}\sum_{\substack{|x|\leq K\sqrt n\\x\equiv n\ (\mathrm{mod}\ 2)}}
 \Pp(S_n=x)H_n(x)\longrightarrow f_a(0)\int_{-K}^Kh_a(z)\dd z.
\]
\end{lemma}
\begin{proof}
Since $\sup_{|x|\leq K\sqrt n}|x|/n^a\leq Kn^{1/2-a}\to0$, the uniform local estimate of \citet[Theorem 3.1]{TV} in Proposition \ref{thm:llt}, together with continuity at zero, gives
\[
 \Pp(S_n=x)=2n^{-a}(f_a(0)+o(1))
\]
uniformly in the sum. Theorem \ref{thm:kernel} gives $H_n(x)=h_a(x/\sqrt n)+o(1)$ uniformly there. The scaled sum equals
\[
 \frac2{\sqrt n}\sum_{\substack{|x|\leq K\sqrt n\\x\equiv n\ (\mathrm{mod}\ 2)}}
 (f_a(0)+o(1))(h_a(x/\sqrt n)+o(1)).
\]
There are at most $K\sqrt n+2$ terms and both errors are uniform. The mesh is $2/\sqrt n$, so the Riemann sum has the asserted limit.
\end{proof}

\begin{lemma}\label{lem:spatial-tail}
If $\Sigma_n(K)$ denotes the sum in \eqref{eq:decomposition} over $|x|>K\sqrt n$, then
\[
 \lim_{K\to\infty}\limsup_{n\to\infty}n^{\delta/2}\Sigma_n(K)=0.
\]
\end{lemma}
\begin{proof}
By \eqref{eq:anti} and Lemma \ref{lem:gaussian},
\[
 n^{\delta/2}\Sigma_n(K)\leq\frac C{\sqrt n}
 \sum_{|x|>K\sqrt n}e^{-cx^2/n}
 \leq C\int_{K/2}^\infty e^{-cu^2}\dd u,\qquad K\geq1.
\]
Integral comparison proves the last inequality after enlarging $C$; the integral tends to zero.
\end{proof}

\begin{proof}[Proof of \eqref{eq:tail} in Theorem \ref{thm:tail}]
Split \eqref{eq:decomposition} at $K\sqrt n$. Lemmas \ref{lem:central} and \ref{lem:spatial-tail}, followed by $K\to\infty$, give
\[
 n^{\delta/2}\Pp(E_n)\longrightarrow f_a(0)\int_\R h_a(z)\dd z.
\]
To evaluate the integral, put $c=\sqrt\delta$ and $\phi(v)=(2\pi)^{-1/2}e^{-v^2/2}$. Tonelli's theorem yields
\begin{align*}
 \int_\R h_a(z)\dd z
 &=\frac4c\int_0^\infty\Phi(-u)\dd u
 =\frac4c\int_0^\infty\int_u^\infty\phi(v)\dd v\dd u\\
 &=\frac4c\int_0^\infty v\phi(v)\dd v
 =\frac4{c\sqrt{2\pi}}=\sqrt{\frac8{\pi\delta}}.
\end{align*}
Almost-sure finiteness and independence of $q$ were proved in Section \ref{sec:model}.
\end{proof}

\begin{corollary}\label{cor:moments}
For $r>0$, $\E[T_0^r]<\infty$ if and only if $r<\delta/2$. In particular $\E[T_0]=\infty$. For $t\geq1$,
\[
 \Pp(T_0/n\geq t\mid E_n)\longrightarrow t^{-\delta/2}.
\]
\end{corollary}
\begin{proof}
For every nonnegative integer-valued $T$, Tonelli's theorem gives
\[
 \E[T^r]=\sum_{n=0}^\infty((n+1)^r-n^r)\Pp(T\geq n+1).
\]
Theorem \ref{thm:tail} makes the summand for $T_0$ comparable to $n^{r-1-\delta/2}$, which is summable exactly when $r<\delta/2$. Also
\[
 \Pp(T_0/n\geq t\mid E_n)=\frac{\Pp(T_0\geq\lceil tn\rceil)}{\Pp(T_0\geq n)}
 \longrightarrow t^{-\delta/2}.
\]
\end{proof}

\section{Position and waiting time given a late return}\label{sec:conditional}
Set $W_n=(\tau_n-n)/n$, with $W_n=\infty$ if $\tau_n=\infty$. The variable $W_n$ measures the elapsed wait in units of $n$, and $1+W_n=\tau_n/n$ is the absolute return time on the same scale. We determine the joint law of position and waiting time, including their dependence. Its density combines the classical Brownian first-passage density \citep[Chapter 2]{MP}, evaluated through the clock $v_a$, with the spatial weights from the LLT of \citet{TV}.

\begin{theorem}\label{thm:joint-pair}
Conditionally on $E_n$, $(S_n/\sqrt n,W_n)\law(Z_a,W_a)$. The limiting pair has joint density
\begin{equation}\label{eq:jointdensity}
 j_a(z,t)=\frac1{\kappa_a}\frac{|z|v_a'(t)}{v_a(t)^{3/2}}
 \phi\left(\frac{|z|}{\sqrt{v_a(t)}}\right),\qquad z\ne0,\ t>0,
\end{equation}
where $v_a'(t)=(1+t)^{-2a}$ and $\phi$ is the standard normal density. Its value at $z=0$ may be set to zero. The spatial marginal has density
\begin{equation}\label{eq:spatialdensity}
 g_a(z)=h_a(z)/\kappa_a.
\end{equation}
The waiting-time marginal has distribution function and density
\begin{align}
 F_a(t)=\Pp(W_a\leq t)&=\sqrt{1-(1+t)^{-\delta}},\qquad t\geq0,\label{eq:waitingcdf}\\
 f_{W_a}(t)&=\frac\delta2\frac{(1+t)^{-2a}}{\sqrt{1-(1+t)^{-\delta}}},\qquad t>0.\label{eq:waitingdensity}
\end{align}
For every $r\geq0$,
\begin{equation}\label{eq:spatialmoments}
 \E\left[\left|S_n/\sqrt n\right|^r\mid E_n\right]
 \longrightarrow \E[|Z_a|^r]
 =\frac{2^{r/2}\Gamma((r+2)/2)}{(r+1)\delta^{r/2}}.
\end{equation}
All these limiting laws are independent of $q$.
\end{theorem}
\begin{proof}
For bounded continuous $\varphi$ and $t>0$, the Markov property gives
\[
 \E[\varphi(S_n/\sqrt n)\ind_{\{W_n\leq t\}}]
 =\sum_{x\equiv n\ (\mathrm{mod}\ 2)}\varphi(x/\sqrt n)\Pp(S_n=x)H_{n,t}(x).
\]
On $|x|\leq K\sqrt n$, apply the LLT \citep[Theorem 3.1]{TV} in \eqref{eq:llt} and Proposition \ref{prop:finite-uniform}; the mesh-$2/\sqrt n$ calculation in Lemma \ref{lem:central} yields
\[
 n^{\delta/2}\sum_{|x|\leq K\sqrt n}\varphi(x/\sqrt n)\Pp(S_n=x)H_{n,t}(x)
 \longrightarrow f_a(0)\int_{-K}^K\varphi(z)h_{a,t}(z)\dd z.
\]
Here and below sums include only the correct parity. Since $H_{n,t}\leq H_n$, the scaled complementary sum is bounded in absolute value by
\[
 C\|\varphi\|_\infty\int_{K/2}^\infty e^{-cu^2}\dd u.
\]
Let $K\to\infty$ and divide by $n^{\delta/2}\Pp(E_n)\to f_a(0)\kappa_a$ to obtain
\begin{equation}\label{eq:subdistribution}
 \E[\varphi(S_n/\sqrt n)\ind_{\{W_n\leq t\}}\mid E_n]
 \longrightarrow\frac1{\kappa_a}\int_\R\varphi(z)h_{a,t}(z)\dd z.
\end{equation}
We next check the time marginal and tightness, so that \eqref{eq:subdistribution} identifies a probability law. The integral calculation in Section \ref{sec:tail}, with $v_a(t)$ in place of $1/\delta$, gives
\[
 \int_\R h_{a,t}(z)\dd z=\sqrt{8v_a(t)/\pi},\qquad
 \frac1{\kappa_a}\int_\R h_{a,t}(z)\dd z=\sqrt{\delta v_a(t)}=F_a(t).
\]
At $t=0$,
\[
 \Pp(W_n=0\mid E_n)=\frac{\Pp(S_n=0)}{\Pp(E_n)}\leq Cn^{-1/2}\longrightarrow0.
\]
Equation \eqref{eq:subdistribution} with $\varphi=1$ gives convergence of the waiting-time distribution functions at every $t>0$. The limit $F_a(t)\downarrow0$ as $t\downarrow0$ excludes an atom at zero in any limiting law. Moreover, $F_a(t)\to1$ as $t\to\infty$ gives tightness of the conditional waiting times, while the spatial tail bound gives tightness of the first coordinate. Equation \eqref{eq:subdistribution}, tested against bounded continuous $\varphi$, determines every subsequential joint limit because its time marginal has no atoms. This proves joint convergence.

For $z\ne0$ and $t>0$,
\[
 \frac\partial{\partial t}h_{a,t}(z)
 =\frac{|z|v_a'(t)}{v_a(t)^{3/2}}
   \phi\left(\frac{|z|}{\sqrt{v_a(t)}}\right).
\]
Thus differentiation of \eqref{eq:subdistribution} gives \eqref{eq:jointdensity}. Letting $t\to\infty$ gives the spatial density $h_a/\kappa_a$; $F_a(0)=0$ excludes mass at time zero, so the joint density accounts for the full law. Differentiating $F_a$ gives \eqref{eq:waitingdensity}.

For the spatial moments, we combine weak convergence with uniform integrability; see \citet[Theorems 3.2.8 and 4.6.3]{Durrett}. The required tail estimate follows from \eqref{eq:anti} and \eqref{eq:gaussian-intro}: for $K\geq1$,
\[
 n^{\delta/2}\E\left[|S_n/\sqrt n|^r
 \ind_{\{|S_n|>K\sqrt n\}}\ind_{E_n}\right]
 \leq C_{a,q,r}\int_{K/2}^\infty u^re^{-cu^2}\dd u.
\]
After division by $n^{\delta/2}\Pp(E_n)$ this tends to zero as $K\to\infty$, uniformly for sufficiently large $n$. Truncation therefore gives convergence of the conditional moments. Finally, Tonelli's theorem and the Gaussian moment integral give, for $c=\sqrt\delta$,
\[
 \int_\R |z|^rh_a(z)\dd z
 =\frac{4c^{-(r+1)}}{r+1}\int_0^\infty v^{r+1}\phi(v)\dd v
 =\frac{2^{(r+3)/2}\Gamma((r+2)/2)}{(r+1)c^{r+1}\sqrt\pi}.
\]
Dividing by $\kappa_a=4/(c\sqrt{2\pi})$ proves \eqref{eq:spatialmoments}. The factor $f_{q,a}(0)=f_a(0)$ cancels in \eqref{eq:subdistribution}.
\end{proof}

\begin{corollary}\label{cor:pair-representation}
For independent $U,R,\varepsilon$ as in the introduction,
\begin{equation}\label{eq:pair-representation}
 (Z_a,W_a)\stackrel{\mathrm d}=
 \left(\frac{\varepsilon UR}{\sqrt\delta},(1-U^2)^{-1/\delta}-1\right).
\end{equation}
For almost every $t>0$, $\E[Z_a^2\mid W_a=t]=2v_a(t)$. For every $z\ne0$, a version of the conditional distribution is
\[
 \Pp(W_a\leq t\mid Z_a=z)=\frac{h_{a,t}(z)}{h_a(z)},\qquad t>0.
\]
\end{corollary}
The change-of-variables verification is given in Appendix~\ref{app:estimates}.

The following moment criterion concerns the limiting waiting time $W_a$; convergence of the conditional moments of $W_n$ would require further estimates on the time tails.

\begin{corollary}\label{cor:waitingmoments}
As $t\to\infty$, $\Pp(W_a>t)\sim\tfrac12t^{-\delta}$. For $r>0$, $\E[W_a^r]<\infty$ if and only if $r<\delta$.
\end{corollary}
\begin{proof}
With $x=(1+t)^{-\delta}$,
\[
 \Pp(W_a>t)=1-\sqrt{1-x}=\frac{x}{1+\sqrt{1-x}}\sim\frac12t^{-\delta}.
\]
The integral $\E[W_a^r]=r\int_0^\infty t^{r-1}\Pp(W_a>t)\dd t$ is finite precisely when $r<\delta$.
\end{proof}

\section{The next return, the last return and escape direction}\label{sec:last}
After the next zero, the walk starts again from the origin at its current absolute time, as in the return-time analysis of \citet[Section 5]{Bertoin}. We first find the limiting law of the last-return time divided by the restart time. Convergence is uniform over late restart times, which yields asymptotic independence of this ratio from the position and the next return. The limit is a negative power of an arcsine variable. Its Brownian antecedent is the last-zero law of \citet{Levy}, and the proof uses the Gaussian-angle calculation presented in \citet[Chapter~5]{MP}.

Under $\Pp_{m,0}$, for even $m\geq2$, define $L_m=\sup\{k\geq m:S_k=0\}$. The event $S_m=0$ has positive probability for the original ERW. Its conditional future law is $\Pp_{m,0}$, so almost-sure escape in Section \ref{sec:model} proves $L_m<\infty$ under this law.

\begin{proposition}\label{prop:restart}
As $m\to\infty$ through even integers, $L_m/m\law V_a$, where
\begin{equation}\label{eq:restart-cdf}
 \Pp(V_a\leq v)=\frac2\pi\arccos(v^{-\delta/2}),\qquad v\geq1.
\end{equation}
Equivalently $V_a\stackrel{\mathrm d}=A^{-1/\delta}$ for an arcsine variable $A$.
\end{proposition}
\begin{proof}
For every integer $k\geq m$, the Markov property gives
\begin{equation}\label{eq:restart-id}
 \Pp_{m,0}(L_m<k)=\E_{m,0}[1-H_k(S_k)].
\end{equation}
Fix $v>1$ and let $k_m=\lfloor mv\rfloor$. Theorem \ref{thm:future-fclt} and Lemma \ref{lem:normalization} imply
\begin{equation}\label{eq:restart-normal}
 \xi_m=\frac{S_{k_m}}{\sqrt{k_m}}\law\sigma_vG,\qquad
 \sigma_v^2=v^{2a-1}v_a(v-1)=\frac{v^\delta-1}\delta,
\end{equation}
where $G$ is standard normal. For $0<K<\infty$,
\begin{align*}
 \E_{m,0}|H_{k_m}(S_{k_m})-h_a(\xi_m)|
 &\leq\sup_{\substack{|x|\leq K\sqrt{k_m}\\(k_m,x)\text{ admissible}}}
       |H_{k_m}(x)-h_a(x/\sqrt{k_m})|\\
 &\quad+\Pp_{m,0}(|\xi_m|>K).
\end{align*}
The first term tends to zero by Theorem \ref{thm:kernel}; tightness from \eqref{eq:restart-normal} removes the second as $K\to\infty$. Thus \eqref{eq:restart-id} tends to
\[
 \E[2\Phi(\sqrt{v^\delta-1}\,|G|)-1].
\]
If $N$ is an independent standard normal and $c=\sqrt{v^\delta-1}$, the Gaussian polar-coordinate argument for the arcsine law \citep[Chapter~5]{MP} gives
\[
 \E[2\Phi(c|G|)-1]=\Pp(|N|\leq c|G|)
 =\frac{4\arctan c}{2\pi}
 =\frac2\pi\arccos(v^{-\delta/2}).
\]
Indeed the event occupies four sectors of angle $\arctan c$ and $1+c^2=v^\delta$. For $1<u<v<w$, eventually
\[
 \{L_m<\lfloor mu\rfloor\}\subset\{L_m/m\leq v\}
 \subset\{L_m<\lfloor mw\rfloor\}.
\]
Let $u\uparrow v$ and $w\downarrow v$. Continuity proves the distributional limit at $v>1$; letting $v\downarrow1$ gives the endpoint, and the expression tends to one at infinity. Finally $\Pp(A\leq u)=(2/\pi)\arcsin\sqrt u$ implies \eqref{eq:restart-cdf} for $A^{-1/\delta}$.
\end{proof}

\begin{proof}[Proof of Theorem~\ref{thm:joint-intro}]
Put $Y_a=1+W_a$. For bounded continuous $g$ on $[1,\infty)$, put $b_m=\E_{m,0}[g(L_m/m)]$ and $b=\E[g(V_a)]$. Proposition \ref{prop:restart} implies
\[
 \eta_n=\sup_{\substack{m\geq n\\m\text{ even}}}|b_m-b|\longrightarrow0.
\]
For bounded continuous $f$ on $\R\times[1,\infty)$, decompose according to the finite value of $\tau_n$ on $E_n$. The deterministic-time Markov property at each such even value gives
\begin{align*}
 &\E[f(S_n/\sqrt n,\tau_n/n)g(T_0/\tau_n)\mid E_n]\\
 &\hspace{12mm}=\E[f(S_n/\sqrt n,\tau_n/n)b_{\tau_n}\mid E_n].
\end{align*}
The difference from $b\E[f(S_n/\sqrt n,\tau_n/n)\mid E_n]$ is bounded directly under the conditional law by $\|f\|_\infty\eta_n$, since $|b_{\tau_n}-b|\leq\eta_n$ on $E_n$. The first two coordinates converge to $(Z_a,Y_a)$ by Theorem \ref{thm:joint-pair}. For the third, the same decomposition gives
\[
 \Pp(T_0/\tau_n>K\mid E_n)
 \leq\sup_{\substack{m\geq n\\m\text{ even}}}\Pp_{m,0}(L_m/m>K).
\]
For fixed $K>1$, Proposition \ref{prop:restart} and continuity at $K$ bound the limsup by $\Pp(V_a>K)$, which vanishes as $K\to\infty$. Hence the triple is tight. Its product test functions have the integrals of the law of $(Z_a,Y_a,V_a)$, with $V_a$ independent of $(Z_a,Y_a)$, which determine every subsequential limit. Corollary~\ref{cor:pair-representation} and Proposition~\ref{prop:restart} identify this law with the representation in \eqref{eq:joint-intro}. Reciprocation is continuous on its support. In particular,
\[
 \Pp(1/V_a\leq r)=\Pp(A\leq r^\delta)=\frac2\pi\arcsin(r^{\delta/2}).
\]
\end{proof}

By the continuous mapping theorem, the joint limit of $(S_n/\sqrt n,\tau_n/n,T_0/n)$ under $E_n$ is consequently
\begin{equation}\label{eq:triple-representation}
 \left(\frac{\varepsilon UR}{\sqrt\delta},
       (1-U^2)^{-1/\delta},\{A(1-U^2)\}^{-1/\delta}\right).
\end{equation}
Its last coordinate has tail $t^{-\delta/2}$ for $t\geq1$, in agreement with Corollary~\ref{cor:moments}.

For example, when $p=7/8$, $a=3/4$ and $\delta=1/2$, the formulas give
\[
 \Pp(E_n)\sim\frac{4f_{3/4}(0)}{\sqrt\pi}n^{-1/4},\qquad
 \E[Z_{3/4}^2]=\frac43.
\]
The medians of $W_{3/4}$ and $V_{3/4}$ are $7/9$ and $4$, respectively. The first measures the waiting time after $n$ relative to $n$; the second measures the last return relative to the next one.

Since $|S_n|\asconv\infty$ and the walk is nearest-neighbor, its eventual sign $D\in\{-1,+1\}$ is almost surely defined by $\sgn(S_n)=D$ for all sufficiently large $n$.

\begin{corollary}\label{cor:direction}
For every $n\geq1$, conditionally on $E_n$, $D$ is symmetric and independent of $(S_n,\tau_n,T_0)$. Moreover,
\begin{align}
 \Pp(S_nD<0)&=\tfrac12\{\Pp(E_n)-\Pp(S_n=0)\},\label{eq:direction-exact}\\
 \Pp(S_nD<0)&\sim f_a(0)\sqrt{\frac2{\pi\delta}}\,n^{-\delta/2}.\label{eq:direction-asymptotic}
\end{align}
\end{corollary}
\begin{proof}
The reflection is applied only after the next zero, leaving the observed past unchanged. Under $\Pp_{m,0}$, reflection of every future position preserves \eqref{eq:transition}, leaves $L_m$ unchanged and exchanges the two signs of $D$. Thus, for Borel $C\subset[m,\infty)$ and $d\in\{-1,+1\}$,
\[
 \Pp_{m,0}(L_m\in C,D=d)=\tfrac12\Pp_{m,0}(L_m\in C).
\]
For bounded measurable $F$, summation over $\{\tau_n=m\}$ gives
\begin{align*}
 &\E[\ind_{E_n}F(S_n,\tau_n,T_0)\ind_{\{D=d\}}]\\
 &=\sum_{\substack{m\geq n\\m\text{ even}}}
 \E\left[\ind_{\{\tau_n=m\}}\E_{m,0}[F(S_n,m,L_m)\ind_{\{D=d\}}]\right]\\
 &=\tfrac12\sum_{\substack{m\geq n\\m\text{ even}}}
 \E\left[\ind_{\{\tau_n=m\}}\E_{m,0}[F(S_n,m,L_m)]\right]
 =\tfrac12\E[\ind_{E_n}F(S_n,\tau_n,T_0)].
\end{align*}
Here $S_n$ is a fixed past parameter in each inner expectation. This proves the conditional independence. If $S_nD<0$, a future zero is required. On $E_n\cap\{S_n\ne0\}$, conditional symmetry makes the opposite-sign probability $1/2$. Since $\{S_n=0\}\subset E_n$, this proves \eqref{eq:direction-exact}. Finally,
\[
 \frac{\Pp(S_n=0)}{\Pp(E_n)}\leq Cn^{-a+\delta/2}=Cn^{-1/2}\longrightarrow0.
\]
Substitution of \eqref{eq:tail} proves \eqref{eq:direction-asymptotic}.
\end{proof}

\section{Last visits to finite targets}\label{sec:finite}
A fixed finite target has vanishing diameter on the $\sqrt n$ scale. For a nearest-neighbor walk started outside its extreme points, entering the set is equivalent to hitting the nearer extreme. This explains why the leading constant does not depend on the number of target sites. To make the argument precise, we first control a fixed target $b$: the martingale boundary is then $bR_{n,m}$, since the transition law is not invariant under spatial translation.

\begin{lemma}\label{lem:fixed-boundary}
Fix $b\in\Z$ and let
\[
 \tau_b^{(n)}=\inf\{m\geq n:S_m=b\},\qquad
 H_n^b(x)=\Pp_{n,x}(\tau_b^{(n)}<\infty).
\]
For every $0<K<\infty$,
\begin{equation}\label{eq:b-kernel}
 \sup_{\substack{|x|\leq K\sqrt n\\(n,x)\text{ admissible}}}
 |H_n^b(x)-h_a(x/\sqrt n)|\longrightarrow0.
\end{equation}
If $B\geq|b|$, there are $c_a,C_a>0$ such that
\begin{equation}\label{eq:b-gaussian}
 H_n^b(x)\leq C_a\exp\left\{-c_a\frac{((|x|-B)_+)^2}{n}\right\}.
\end{equation}
\end{lemma}
\begin{proof}
Let $x_n/\sqrt n\to z\in\R$ and fix $T>0$. Consider
\[
 Z_n^b(t)=\frac{R_{n,n+\lfloor nt\rfloor}(S_{n+\lfloor nt\rfloor}-b)}{\sqrt n}
 =\frac{M_{n+\lfloor nt\rfloor}^{(n)}}{\sqrt n}
   -\frac{bR_{n,n+\lfloor nt\rfloor}}{\sqrt n}.
\]
The second term is bounded uniformly by $|b|/\sqrt n$, so Theorem \ref{thm:future-fclt} gives $Z_n^b\law z+B(v_a(\cdot))$. Along subsequences with $x_n\geq b$, a visit to $b$ is exactly the event $\inf_{t\leq T}Z_n^b(t)\leq0$; along those with $x_n\leq b$, use $\sup_{t\leq T}Z_n^b(t)\geq0$. The nearest-neighbor property justifies both identities. If $z\ne0$, the initial side is eventually determined by its sign. If $z=0$, both extrema arguments give probability one, since the Brownian minimum and maximum have no atom at zero. Thus the proof of Lemma \ref{lem:finite-hitting} yields, for every $z$,
\begin{equation}\label{eq:b-finite}
 \Pp_{n,x_n}(\tau_b^{(n)}\leq n+\lfloor nT\rfloor)\longrightarrow h_{a,T}(z).
\end{equation}
The limit thus follows from the vanishing normalized boundary $bR_{n,m}/\sqrt n$.

For the infinite horizon put $m_T=n+\lfloor nT\rfloor$. On $\{m_T<\tau_b^{(n)}<\infty\}$ either
\[
 |M_{m_T}^{(n)}|\leq |b|+\eta\sqrt n
 \quad\text{or}\quad
 \sup_{m\geq m_T}|M_m^{(n)}-M_{m_T}^{(n)}|\geq\eta\sqrt n,
\]
because $|M_{\tau_b^{(n)}}^{(n)}|=|b|R_{n,\tau_b^{(n)}}\leq|b|$. The second probability is bounded by \eqref{eq:Doob-tail}. In the first, $|b|/\sqrt n\to0$, and the Gaussian interval argument in \eqref{eq:small-interval} gives a uniform limsup bound $C_a\eta$. Therefore the late-hit probability has the same bound $C_a\{\eta+\eta^{-2}(1+T)^{-\delta}\}$ as before. Let $T\to\infty$ with $\eta=(1+T)^{-\delta/4}$. Combined with \eqref{eq:b-finite} and a compactness argument, this proves \eqref{eq:b-kernel}.

Finally a hit requires
\[
 |M_{\tau_b^{(n)}}^{(n)}-M_n^{(n)}|
 =|bR_{n,\tau_b^{(n)}}-x|\geq(|x|-|b|)_+.
\]
The maximal inequality \eqref{eq:Azuma} and the bound on the sum of squared increments prove \eqref{eq:b-gaussian}.
\end{proof}

\begin{proof}[Proof of the finite-set assertions in Theorem \ref{thm:tail}]
Let $b_- =\min F$, $b_+=\max F$, $B=\max_{b\in F}|b|$, and
\[
 H_n^F(x)=\Pp_{n,x}(S_m\in F\text{ for some }m\geq n).
\]
If $x>b_+$, the first point of $F$ that can be hit is $b_+$, so $H_n^F(x)=H_n^{b_+}(x)$. If $x<b_-$, similarly $H_n^F(x)=H_n^{b_-}(x)$. For $b_-\leq x\leq b_+$, choose any fixed $b_0\in F$; then
\[
 H_n^{b_0}(x)\leq H_n^F(x)\leq1.
\]
Lemma \ref{lem:fixed-boundary} gives $H_n^{b_0}(x)\to1$ uniformly over these finitely many states of the appropriate parity, while $h_a(x/\sqrt n)\to1$. Consequently
\[
 \sup_{\substack{|x|\leq K\sqrt n\\(n,x)\text{ admissible}}}
 |H_n^F(x)-h_a(x/\sqrt n)|\longrightarrow0.
\]
A hit of $F$ requires a martingale displacement of at least $(|x|-B)_+$, so
\[
 H_n^F(x)\leq C_a\exp\{-c_a((|x|-B)_+)^2/n\}.
\]
For $|x|>K\sqrt n$ and sufficiently large $n$, with $K\geq1$ fixed, $|x|-B\geq|x|/2$. Hence the same integrable Gaussian tail used in Lemma \ref{lem:spatial-tail} applies. Since, for $n\geq1$,
\[
 \Pp(T_F\geq n)=\sum_{x\equiv n\ (\mathrm{mod}\ 2)}\Pp(S_n=x)H_n^F(x),
\]
the central-window Riemann sum and the spatial-tail estimate give
\[
 n^{\delta/2}\Pp(T_F\geq n)\longrightarrow f_a(0)\int_\R h_a(z)\dd z=f_a(0)\kappa_a.
\]
Finiteness follows from $|S_n|\asconv\infty$. Applying the integer moment-sum formula in Corollary \ref{cor:moments} proves the asserted criterion for $T_F$.
\end{proof}

The argument is uniform for targets whose location and diameter are negligible on the $\sqrt n$ scale. The target may depend on the observation time $n$, but is held fixed throughout the future defining $T_F$.

\begin{corollary}\label{cor:subdiffusive-targets}
Let $(r_n)$ be a deterministic nonnegative sequence with $r_n=o(\sqrt n)$. Then
\begin{equation}\label{eq:subdiffusive-targets}
 \sup_{\varnothing\ne F\subseteq[-r_n,r_n]\cap\Z}
 \left|n^{a-1/2}\Pp(T_F\geq n)-f_a(0)\kappa_a\right|
 \longrightarrow0.
\end{equation}
The supremum is over deterministic sets. No bound on the cardinality of $F$, other than its inclusion in the stated interval, is required.
\end{corollary}
\begin{proof}
Fix $0<K<\infty$ and put $m_T=n+\lfloor nT\rfloor$. First extend Lemma \ref{lem:fixed-boundary} uniformly to $|b|\leq r_n$. For any such $b=b_n$,
\[
 \sup_{m\geq n}\frac{|bR_{n,m}|}{\sqrt n}\leq\frac{r_n}{\sqrt n}\longrightarrow0.
\]
For every sequence $x_n/\sqrt n\to z$, the finite-horizon minimum or maximum argument in that lemma therefore gives the limit $h_{a,T}(z)$. For the late-hit estimate, its deterministic inclusion now gives
\[
 \Pp_{n,x}(m_T<\tau_b^{(n)}<\infty)
 \leq\Pp_{n,x}(|M_{m_T}^{(n)}|\leq r_n+\eta\sqrt n)
       +C_a\eta^{-2}(1+T)^{-\delta}.
\]
For fixed $T\geq1$ and $\eta>0$, the limsup of the first term, uniformly over $|x|\leq K\sqrt n$, is at most $C_a\eta$: use $r_n/\sqrt n\to0$ and the normal-density argument of \eqref{eq:small-interval}. Choosing $\eta=(1+T)^{-\delta/4}$ and then letting $T\to\infty$ proves the infinite-horizon limit. Compactness of the scaled starting points gives
\[
 \sup_{\substack{|b|\leq r_n,\ b\in\Z\\|x|\leq K\sqrt n,\ (n,x)\text{ admissible}}}
 |H_n^b(x)-h_a(x/\sqrt n)|\longrightarrow0.
\]
For any nonempty $F\subseteq[-r_n,r_n]\cap\Z$, outside its extreme points $H_n^F$ equals the kernel of the nearer extreme. Between them, choose $b\in F$ and use $H_n^b(x)\leq H_n^F(x)\leq1$; there $|x|/\sqrt n\leq r_n/\sqrt n\to0$. Thus the same convergence holds uniformly over both $F$ and the central window. A hit of $F$ also requires a martingale displacement of at least $(|x|-r_n)_+$, so \eqref{eq:Azuma} gives
\[
 H_n^F(x)\leq C_a\exp\{-c_a((|x|-r_n)_+)^2/n\}.
\]
For $|x|>K\sqrt n$, $K\geq1$ fixed, this is bounded by $C_a e^{-c_a x^2/(4n)}$ for all sufficiently large $n$. The central Riemann sum and the integrable spatial bound in Lemmas \ref{lem:central} and \ref{lem:spatial-tail} are consequently uniform in $F$, and prove \eqref{eq:subdiffusive-targets}.
\end{proof}

\section{Mean local times and occupation}\label{sec:occupation}
The following occupation estimates are consequences of the uniform LLT of \citet[Theorem 3.1]{TV}, recorded in Proposition~\ref{thm:llt}. Almost-sure finiteness follows from the known transience, as recalled in Section~\ref{sec:model}. These estimates quantify the divergence of the mean occupation although the occupation is almost surely finite.
For $x\in\Z$, define
\[
 \ell_N(x)=\sum_{n=0}^N\ind_{\{S_n=x\}},\qquad
 \ell_\infty(x)=\lim_{N\to\infty}\ell_N(x).
\]

\begin{corollary}\label{thm:localtime}
For every fixed $x\in\Z$,
\[
 \E[\ell_N(x)]\sim\frac{f_a(0)}{1-a}N^{1-a}
 =\frac{f_a(0)}{2-2p}N^{2-2p}.
\]
Moreover $\ell_\infty(x)<\infty$ almost surely and $\E[\ell_\infty(x)]=\infty$. At the origin this remains true after excluding the initial visit, so the total number of returns $N_0=\sum_{n=1}^{\infty}\ind_{\{S_n=0\}}$ is finite almost surely with infinite expectation.
\end{corollary}
\begin{proof}
At times $n\equiv x\pmod2$, the LLT of \citet[Theorem 3.1]{TV}, in the form \eqref{eq:llt}, gives
\[
 \Pp(S_n=x)=2n^{-a}(f_a(0)+o(1)),
\]
since $x/n^a\to0$; at the other times the probability is zero. Thus
\[
 \E[\ell_N(x)]\sim2f_a(0)\sum_{\substack{1\leq n\leq N\\n\equiv x\ (\mathrm{mod}\ 2)}}n^{-a}
 \sim\frac{f_a(0)}{1-a}N^{1-a}.
\]
The factor $2$ in the local limit theorem is canceled by summing over one parity class. Almost-sure escape proves finiteness of $\ell_\infty(x)$. Monotone convergence and the divergent mean asymptotic prove its infinite expectation.
\end{proof}

The same LLT \citep[Theorem 3.1]{TV}, now summed at sites on the superdiffusive scale, gives the following mean occupation profile.\nopagebreak[4]

\begin{proposition}\label{thm:occupation}
Let $x_N\in\Z$ satisfy $x_N/N^a\to y\in\R$. Then
\begin{equation}\label{eq:occupation}
 N^{a-1}\E[\ell_N(x_N)]\longrightarrow\rho_{q,a}(y),\qquad
 \rho_{q,a}(y)=\int_0^1t^{-a}f_{q,a}(y/t^a)\dd t.
\end{equation}
The function $\rho_{q,a}$ is nonnegative and has integral one.
\end{proposition}
\begin{proof}
For $x\equiv k\pmod2$, write the uniform local estimate \eqref{eq:llt} from \citet[Theorem 3.1]{TV} as
\[
 \Pp(S_k=x)=2k^{-a}\{f_{q,a}(x/k^a)+\epsilon_k(x)\},\qquad
 \sup_{x\equiv k\ (\mathrm{mod}\ 2)}|\epsilon_k(x)|\longrightarrow0.
\]
Fix $\eta\in(0,1)$. For $\eta N\leq k\leq N$, the arguments $x_N/k^a$ remain in a compact set. Uniform continuity there and $x_N/N^a\to y$ show that
\begin{align*}
 &N^{a-1}\sum_{\substack{\eta N\leq k\leq N\\k\equiv x_N\ (\mathrm{mod}\ 2)}}\Pp(S_k=x_N)\\
 &\quad=\frac2N\sum_{\substack{\eta N\leq k\leq N\\k\equiv x_N\ (\mathrm{mod}\ 2)}}
 (k/N)^{-a}f_{q,a}\left(\frac{x_N/N^a}{(k/N)^a}\right)+o(1)
 \longrightarrow\int_\eta^1t^{-a}f_{q,a}(y/t^a)\dd t.
\end{align*}
The accumulated local-limit error is at most $C_\eta\sup_{k\geq\eta N,x}|\epsilon_k(x)|\to0$. The mesh is $2/N$, which cancels the span factor $2$. For the omitted times \eqref{eq:anti} gives
\[
 N^{a-1}\sum_{1\leq k<\eta N}\Pp(S_k=x_N)
 \leq CN^{a-1}\sum_{1\leq k<\eta N}k^{-a}\leq C\eta^{1-a}.
\]
The term at $k=0$ is at most $N^{a-1}\to0$. The limiting integrand is bounded by $\|f_{q,a}\|_\infty t^{-a}$, integrable at zero. Let $\eta\downarrow0$ to prove \eqref{eq:occupation}. Finally, Tonelli's theorem and $u=y/t^a$ give
\[
 \int_\R\rho_{q,a}(y)\dd y
 =\int_0^1\int_\R f_{q,a}(u)\dd u\dd t=1.
\]
\end{proof}

\section{Further questions}
The distribution of the total number of visits to a fixed finite set remains a separate question: its infinite expectation does not determine its tail. Further problems include quantitative error bounds in the last-return asymptotic and targets whose location or diameter is comparable to $\sqrt n$, beyond the range of Corollary~\ref{cor:subdiffusive-targets}. Time-dependent targets, whose sites change during the trajectory, also require a different hitting analysis. The present estimates keep $a\in(1/2,1)$ fixed and do not address the transition as $a\downarrow1/2$ or $a\uparrow1$. Extending this approach to other memory rules would require suitable local estimates and control of future hitting probabilities.

\clearpage
\appendix
\section{Auxiliary estimates}\label{app:estimates}
The proofs collected here supply the product comparison, the second-moment bound and the change of variables used in the main argument.

\begin{proof}[Proof of Lemma~\ref{lem:normalization}]
The inequalities $0\leq u-\log(1+u)\leq u^2/2$ for $u\geq0$ and integral comparison give
\[
 \log R_{n,m}=-a\sum_{j=n}^{m-1}\frac1j+e_{n,m},\qquad
 |e_{n,m}|\leq\frac{a^2}2\sum_{j=n}^\infty j^{-2},\qquad
 0\leq\sum_{j=n}^{m-1}\frac1j-\log\frac mn\leq\frac1n.
\]
In particular $|\log R_{n,m}+a\log(m/n)|\leq C_a/n$. Exponentiating gives the two-sided estimates. Since $\sup_{t\leq T}|(n+\lfloor nt\rfloor)/n-(1+t)|\leq1/n$, the same estimates prove \eqref{eq:Rlimit}. If $M=n+\lfloor nT\rfloor$, then
\[
 \frac1n\sum_{m\geq M}R_{n,m+1}^2
 \leq C_a n^{2a-1}\sum_{m\geq M}(m+1)^{-2a}
 \leq\frac{C_a n^{2a-1}M^{1-2a}}{2a-1}
 \leq C_a(1+T)^{1-2a}.
\]
Here $M\geq n(1+T)/2$ and $2a>1$ were used. The exact product formula is $R_{n,m}=\Gamma(m)\Gamma(n+a)/\{\Gamma(m+a)\Gamma(n)\}$.
\end{proof}

\begin{proof}[Proof of Lemma~\ref{lem:second}]
The conditional second-moment identity in \citet[(A.3)]{Bercu} follows by expanding $(S_m+X_{m+1})^2$; taking expectations recovers the recursion already derived by \citet[(14)]{ST}:
\[
 \E[S_{m+1}^2\mid\F_m]=(1+2a/m)S_m^2+1.
\]
Thus $u_m=\E_{n,x}[S_m^2]$ satisfies
\[
 u_m=x^2\prod_{j=n}^{m-1}(1+2a/j)
       +\sum_{k=n}^{m-1}\prod_{j=k+1}^{m-1}(1+2a/j).
\]
The same logarithmic calculation as in the proof of Lemma \ref{lem:normalization}, with $2a$ in place of $a$, bounds each product by $C_a(m/(k+1))^{2a}$, and the first by $C_a(m/n)^{2a}$. On the stated interval these ratios are at most $T+2$. Hence
\[
 u_m\leq C_{a,T}K^2n+C_{a,T}(nT+1)\leq Cn.
\]
\end{proof}

\begin{proof}[Proof of Corollary~\ref{cor:pair-representation}]
Set $u(t)=\sqrt{\delta v_a(t)}$, so that $u'(t)=\sqrt\delta\,v_a'(t)/(2\sqrt{v_a(t)})$. Given $U=u$, the density of $Z'=\varepsilon uR/\sqrt\delta$ is
\[
 f_{Z'\mid U=u}(z)=\frac{\delta|z|}{2u^2}
 \exp\left(-\frac{\delta z^2}{2u^2}\right).
\]
Changing variables from $u$ to $t$ gives the joint density
\[
 \frac{|z|}{2v_a(t)}e^{-z^2/(2v_a(t))}
 \frac{\sqrt\delta\,v_a'(t)}{2\sqrt{v_a(t)}}
 =\frac{\sqrt\delta\,|z|v_a'(t)}{4v_a(t)^{3/2}}e^{-z^2/(2v_a(t))}
 =j_a(z,t),
\]
since $1/(\kappa_a\sqrt{2\pi})=\sqrt\delta/4$. This proves equality of the joint laws. Given the second coordinate $t$, $U^2=\delta v_a(t)$ and $R$ remains independent. Since $\E[R^2]=2$, the second-moment formula follows. Dividing the subdistribution in \eqref{eq:subdistribution} by the spatial density proves the last assertion.
\end{proof}

\end{document}